\documentclass{article}
\usepackage{amsfonts}
\usepackage{amsmath}
\usepackage{ulem}
\usepackage{color}

\newtheorem{theorem}{Theorem}

\newtheorem{lemma}[theorem]{Lemma}

\begin{document}

\title{Critical branching processes in random environment with immigration: the size of the only surviving family\thanks{This work is supported by the Russian Science
Foundation under the grant 19-11-00111}}
\author{C. Smadi\thanks{%
Department of Discrete Mathematics, Steklov Mathematical Institute of
Russian Academy of Sciences, 8 Gubkin Street, 117 966 Moscow GSP-1, Russia
E-mail: charline.smadi@inrae.fr},\thinspace\ V. A. Vatutin\thanks{%
Department of Discrete Mathematics, Steklov Mathematical Institute of
Russian Academy of Sciences, 8 Gubkin Street, 117 966 Moscow GSP-1, Russia,
E-mail: vatutin@mi-ras.ru}}
\date{\today}
\maketitle

\begin{abstract}
We consider a critical branching process $Y_{n}$ in an i.i.d. random
environment, in which one immigrant arrives at each generation. Let $%
\mathcal{A}_{i}(n)$ be the event that all individuals alive at time $n$ are
offspring of the immigrant which joined the population at time $i$. We study
the conditional distribution of $Y_{n}$ given $\mathcal{A}_{i}(n)$ when $n$
is large and $i$ follows different asymptotics which may be related to $n$ ($%
i$ fixed, close to $n$, or going to infinity but far from $n$).


\noindent \textbf{AMS 2000 subject classifications.} Primary 60J80;
Secondary 60G50.\newline

\noindent \textbf{Keywords.} Branching process, random environment,
immigration, conditioned random walk
\end{abstract}

\section{Introduction and main result}

We consider a branching process with immigration evolving in a random
environment. Individuals in this process reproduce
independently of each other according to random offspring distributions
which vary from one generation to the other. In addition, an immigrant
enters the population at each generation. A formal definition of such a
process looks as follows. Let $\Delta $ be the space of all probability
measures on $\mathbf{N}_{0}:=\{0,1,2,\ldots \}.$ Equipped with a metric, $\Delta $ is a Polish space. Let $F$ be a random variable taking
values in $\Delta $, and let $F_{n},n\in \mathbf{N}:=\mathbf{N}%
_{0}\backslash \left\{ 0\right\} $ be a sequence of independent copies of $F$%
. The infinite sequence $\mathcal{E}=\left\{ F_{n},n\in \mathbf{N}\right\} $
is called a random environment.

A sequence of $\mathbf{N}_{0}$-valued random variables $\mathbf{Y}=\left\{
Y_{n},\ n\in \mathbf{N}_{0}\right\} $ specified on\ the respective
probability space $(\Omega ,\mathcal{F},\mathbf{P})$ is called a branching
process with one immigrant joining each generation and evolving in  random environment (BPIRE for short), if $Y_{0}=1$ and,
given $\mathcal{E}$ the process $\mathbf{Y}$ is a Markov chain with
\begin{equation*}
\mathcal{L}\left( Y_{n}|Y_{n-1}=y_{n-1},F_{i}=f_{i},i=1,2,\ldots \right) =%
\mathcal{L}(\xi _{n1}+\ldots +\xi _{ny_{n-1}}+1)  \label{BasicDefBPimmigr}
\end{equation*}%
for every $n\in \mathbf{N}$, $y_{n-1}\in \mathbf{N}_{0}$ and $%
f_{1},f_{2},...\in \Delta $, where $\xi _{n1},\xi _{n2},\ldots $ are i.i.d.
random variables with distribution $f_{n}.$
We assume in the sequel that if $Y_{n-1}=y_{n-1}>0$ is the population size
of the ($n-1)$th generation of $\mathbf{Y}$ then first $\xi _{n1}+\ldots
+\xi _{ny_{n-1}}$ individuals of the $n$th generation are born and
afterwards one immigrant enters the population.\newline

The tail distribution of the life-periods of BPIRE's were considered in \cite{DLVZ19} and \cite{LVZ19} under weaker assumptions than the ones we impose here.
An $(i,n)$-clan of a BPIRE is the set of individuals alive at generation $n$
and being descendants of the immigrant which entered the population at
generation $i$. We say that only the $(i,n)$-clan survives in $\mathbf{Y}$
at moment $n$ if $Y_{n}^{-}:=\xi _{n1}+\ldots +\xi _{ny_{n-1}}>0$ and all
the $Y_{n}^{-}$ particles belong to the $(i,n)$-clan.

Let $\mathcal{A}_{i}(n)$ be the event that only the $(i,n)$-clan survives in
$\mathbf{Y}$ at moment $n$. The asymptotic behavior of the probability $%
\mathbf{P}\left( \mathcal{A}_{i}(n)\right) $ as $n\rightarrow \infty $ and $%
i $ varies with $n$ in an appropriate way was investigated in \cite{SV2021} for the critical BPIRE's (see Assumption \textbf{A2}) and in \cite{VD2020} for the subcritical BPIRE's .
The present paper complements the results of \cite{SV2021} by describing
the distribution of the population size of the process at moment $n$ given
the event $\mathcal{A}_{i}(n)$. {The estimation of the population size in this case} is important from an evolutionary point of view because it provides information on the probability of survival of the population in a variable environment. Successive migrants bring genetic diversity to the population, and when there is only one surviving lineage, this implies that the population is poor in genetic diversity. If the population size of the process at moment $n$, given
the event $\mathcal{A}_{i}(n)$, is small, the population will therefore be very vulnerable to an environmental change that could make the present genetic type less adapted.
The question of quantifying the size and the genetic diversity of a population is also related to the concept of founder effect in population genetics.
It is the loss of genetic variation that occurs when a new population is established by a very small number of individuals from a larger population. Because of the loss of genetic variation, the new population may be distinctly different, both genotypically and phenotypically, from the parent population from which it was derived. In extreme cases, the founder effect is thought to lead to speciation and subsequent evolution of new species (see \cite{matzke2014model} for more details on founder effect).

We need to consider, along with the process $\mathbf{Y}$, a standard
branching process $\mathbf{Z}=\left\{ Z_{n},\ n\in \mathbf{N}_{0}\right\} $
in the random environment $\mathcal{E}$ which, given $\mathcal{E}$ is a
Markov chain with $Z_{0}=1$ and
\begin{equation*}
\mathcal{L}\left( Z_{n}|Z_{n-1}=z_{n-1},F_{i}=f_{i},i=1,2,\ldots \right) =%
\mathcal{L}(\xi _{n1}+\ldots +\xi _{nz_{n-1}})  \label{BPordinary}
\end{equation*}%
for $n\in \mathbf{N}$, $z_{n-1}\in \mathbf{N}_{0}$ and $f_{1},f_{2},...\in
\Delta $.

To formulate our results we introduce the so-called associated random walk $%
\mathbf{S}=\left\{ S_{n},n\in \mathbf{N}_{0}\right\} $ (see \cite{4h} for instance). This random walk has
increments $X_{n}=S_{n}-S_{n-1}$, $n\geq 1$, defined as
\begin{equation*}
X_{n}=\log m\left( F_{n}\right)
\end{equation*}%
which are i.i.d. copies of the logarithmic mean offspring number $X:=\log $ $%
m(F)$ with%
\begin{equation*}
m(F):=\sum_{j=1}^{\infty }jF\left( \left\{ j\right\} \right) .
\end{equation*}

With each measure $F$ we associate the respective probability generating
function
\begin{equation*}
F(s):=\sum_{j=0}^{\infty }F\left( \left\{ j\right\} \right) s^{j}.
\end{equation*}%
Introduce the following assumptions:\\

\noindent \textbf{Hypothesis A1}. The probability generating function $F(s)$
is geometric with probability 1:%
\begin{equation}
F(s)=\frac{q}{1-ps}=\frac{1}{1+m(F)(1-s)}  \label{Frac_generating}
\end{equation}%
with random $p,q\in (0,1)$ satisfying $p+q=1,$ and the random variable
\begin{equation*}
X:=\log m(F)=\log \frac{p}{q}
\end{equation*}%
has a nonlattice distribution.\\

\noindent \textbf{Hypothesis A2}. The branching process $\mathbf{Z}$ is
critical: $\mathbf{E}\left[ X\right] =0,$ and
\begin{equation*}
\mathbf{E}\left[ e^{X}+e^{-X}\right] <\infty .
\end{equation*}

\noindent \textbf{Hypothesis A3}. The distribution of $X$ is continuous.%
\newline

Denote $Z_{i,n}$ the size of the $(i,n)$-clan at moment $n$ and, introducing the associated random walk $S$ with $S_0=0,$ set
\begin{equation*}
Y_{i,n}:=e^{S_{i}-S_{n}}Z_{i,n}.
\end{equation*}

The main result of the present paper is the following theorem.

\begin{theorem}
\label{T_conditional}If Hypotheses A1--A2 are valid then

1) for any fixed $N$ and $s\in \lbrack 0,1]$
\begin{equation}
\lim_{n\rightarrow \infty }\mathbf{E}\left[ s^{Z_{n-N,n}}|\mathcal{A}%
_{n-N}(n)\right] =:\Theta _{N}(s)=\mathbf{E}\left[ s^{\vartheta _{N}}\right]
,  \label{YagDiscrete}
\end{equation}%
where $\vartheta _{N}$ is a proper nondegenerate random variable;

2) for any fixed $i$ and $\mathcal{\beta }\geq 0$
\begin{equation*}
\lim_{n\rightarrow \infty }\mathbf{E}\left[ e^{-\mathcal{\beta }Y_{i,n}}|\mathcal{A}_{i}(n)\right] =:\Lambda _{i}\left( \mathcal{\beta }%
\right) =\mathbf{E}\left[ e^{-\beta \hat{Y}_{i,\infty }}\right] ,
\label{YagSmall}
\end{equation*}%
where $\hat{Y}_{i,\infty }$ a proper strictly positive random variable;

3) if $\min \left( i,n-i\right) \rightarrow \infty $ and, in addition,
Hypotheses A3 is valid then for any $\mathcal{\beta }\geq 0$
\begin{equation*}
\lim_{n\rightarrow \infty }\mathbf{E}\left[ e^{-\mathcal{\beta }Y_{i,n}}|\mathcal{A}_{i}(n)\right] =:\Lambda \left( \mathcal{\beta }%
\right) =\mathbf{E}\left[ e^{-\beta \hat{Y}}\right] ,  \label{YagInterm}
\end{equation*}%
where $\hat{Y}$ a proper strictly positive random variable.
\end{theorem}

Roughly speaking, Theorem \ref{T_conditional} establishes that on the event $\mathcal{A}_{i}(n)$, the population size $Z_{i,n}$ behaves as $e^{S_n-S_i}$.
Previous works have shown that, conditioned on the event of survival at time $n$ of a  Galton-Watson process {evolving in an environment with independent identically distribute components}, the population size at this moment behaves as $e^{S_n-S_{\tau(n)}}$, where $S_{\tau(n)}$ is the minimum of the random walk $S$ on $\{0,...,n\}$ (see for instance Theorem 1.3 in \cite{4h}, or Theorem 1.4 in \cite{ABKV} or \cite{ABKV2}).
As previously observed under different assumptions on the random environment (see, for instance, \cite{VDS2013} for a comprehensive review on the critical and subcritical cases (before 2013) or the recent monograph \cite{GV2017}) the survival of a branching process in random environment until a {distant} time $n$ is essentially determined by its survival until the moment when the associated random walk $\mathbf{S}$ attains its infimum on the interval $[0,n]$.
{The idea is that if we divide the trajectory of the process on the interval $[0,n]$ into two parts,  before the running infimum $\tau(n)$ of the associated random walk $\mathbf{S}$ on $[0,n]$, and after this moment, the process will live in a favorable environment after the moment $\tau(n)$, and will thus survive with a non-negligible probability until time $n$, provided it survived until $\tau(n)$.}
{For the process with immigration we consider here, moment $\tau(n)$ should be, as a rule, close to close time $i$.} Indeed, on one hand, on the event $\mathcal{A}_{i}(n)$, every family {generated by an immigrant joined the population} before time $i$ is extinct at the time $n$, and thus has undergone bad environments before {observation} time $n$. On the other hand, the fact that the $(i,n)$-clan is nonempty at time $n$ implies that $\inf \{S_k-S_i, i \leq k \leq n\}$ is likely to be non-negative. Finally, the difference $S_n-S_i$ is also likely to be not too big, otherwise some of the immigrants arrived after time $i$ would have a positive line of descents at time $n$. These {nonrigorous} considerations indicate that on the event $\mathcal{A}_{i}(n)$:
\begin{itemize}
\item $S_n-S_i$ is likely to be of the same order  {as} $S_n-S_{\tau(n)}$
\item $S_n-S_i$ is likely to be relatively small
\end{itemize}
 {The results of Theorem \ref{T_conditional} confirm to a certain extend these hypothesis and similar in spirit to Theorem 1.3 in \cite{4h},  Theorem 1.4 in \cite{ABKV} and paper \cite{ABKV2}.}

The rest of the paper is organised as follows. In Section \ref{sec_aux_res}
we collect some auxiliary results dealing with the probability of the event $%
\mathcal{A}_{i}(n)$. Section \ref{sec_nN} is dedicated to the proof of point
1) of Theorem \ref{T_conditional}. The proof of point 3) of Theorem \ref%
{T_conditional} is provided in Section \ref{sec_ini_to_infty}. Finally, the
proof of Theorem~\ref{T_conditional} is completed in Section \ref%
{sec_fixed_i} by considering the case when $i$ is fixed.\newline

In the sequel we will denote by $C,C_{1}, C_{2},...$ constants which may
vary from line to line and by $K_{1},K_{2},...$ some fixed constants.

\section{Auxiliary results}

\label{sec_aux_res}

Given the environment $\mathcal{E}=\left\{ F_{n},n\in \mathbf{N}\right\} $,
we {introduce the i.i.d. sequence of generating functions
\begin{equation*}
F_{n}(s):=\sum_{j=0}^{\infty }F_{n}\left( \left\{ j\right\} \right)
s^{j},\quad s\in \lbrack 0,1],
\end{equation*}%
and use below the cmpositions of $F_{1},...,F_{n}$ specified for $0\leq
i\leq n-1$ by the equalities
\begin{equation*}
\left\{
\begin{array}{l}
F_{i,n}(s):=F_{i+1}(F_{i+2}(\ldots F_{n}(s)\ldots )),\quad \\
F_{n,i}(s):=F_{n}(F_{n-1}(\ldots F_{i+1}(s)\ldots )),%
\end{array}%
\right.
\end{equation*}%
and $F_{n,n}(s):=s$ for $i=n$. }

Set
\begin{equation*}
\mathcal{H}_{i,n}(s):=\left( 1-F_{i,n}(s)\right) \prod_{j\neq
i}^{n-1}F_{j,n}(0),\quad \mathcal{H}_{i,n}:=\mathcal{H}_{i,n}(0).
\end{equation*}

For $1\leq i\leq n$ introduce the notation
\begin{eqnarray}
&&a_{i,n}:= e^{S_{i}-S_{n}},\,\,a_{n}:=a_{0,n}=e^{-S_{n}},\,b_{0}:=0,
\label{defainbin} \\
&&b_{i,n}:=\sum_{k=i}^{n-1}e^{S_{i}-S_{k}},\
b_{n}:=b_{0,n}=\sum_{k=0}^{n-1}e^{-S_{k}}=:1+B_{1,n}.  \notag
\end{eqnarray}

We can check by induction that if condition (\ref{Frac_generating}) is valid
then for $i\leq n$

\begin{equation}
1-F_{i,n}(s)=\frac{a_{i}}{a_{n}\left( 1-s\right) ^{-1}+b_{n}-b_{i}}.
\label{expr_Fin1}
\end{equation}%
In particular,
\begin{equation}
F_{i,n}(0)=\frac{a_{n}+b_{n}-b_{i+1}}{a_{n}+b_{n}-b_{i}}.  \label{Frac_F}
\end{equation}

Now we provide an expression in terms of $a_{k}$ and $b_{k},\, k\geq 0,$ for the
random variable $\mathcal{H}_{i,n}(s).$

\begin{lemma}
\label{L_fractional} Under Hypothesis $A1$ for any $i=0,1,...,n-1$
\begin{eqnarray*}
\mathcal{H}_{i,n}(s) &=&\frac{1}{a_{i,n}(1-s)^{-1}+b_{i,n}}\frac{%
a_{i,n}+b_{i,n}}{a_{i,n}+b_{i,n}-1}\frac{a_{n}}{a_{n}+b_{n}} \\
&=&\frac{a_{i}}{a_{n}(1-s)^{-1}+b_{n}-b_{i}}\frac{a_{n}+b_{n}-b_{i}}{%
a_{n}+b_{n}-b_{i+1}}\frac{a_{n}}{b_{n+1}}.
\end{eqnarray*}%
In particular,
\begin{equation*}
\mathbf{P}\left( \mathcal{A}_{i}(n)|\mathbf{S}\right) =\mathbf{E}\left[
\mathcal{H}_{i,n}(0)\Bigg|\mathbf{S}\right] =\frac{a_{i}}{a_{n}+b_{n}-b_{i+1}%
}\frac{a_{n}}{b_{n+1}}.  \label{ProbabA_i}
\end{equation*}
\end{lemma}

\textbf{Proof}. The desired statements are direct consequences of (\ref%
{expr_Fin1}) and (\ref{Frac_F}):
\begin{eqnarray*}
\mathcal{H}_{i,n}(s) &=&\frac{\left( 1-F_{i,n}(s)\right) }{F_{i,n}(0)}%
\prod_{j=0}^{n-1}F_{j,n}(0) \\
&=&\frac{a_{i}}{a_{n}(1-s)^{-1}+b_{n}-b_{i}}\frac{a_{n}+b_{n}-b_{i}}{%
a_{n}+b_{n}-b_{i+1}}\prod_{j=0}^{n-1}\frac{a_{n}+b_{n}-b_{j+1}}{%
a_{n}+b_{n}-b_{j}} \\
&=&\frac{a_{i}}{a_{n}(1-s)^{-1}+b_{n}-b_{i}}\frac{a_{n}+b_{n}-b_{i}}{%
a_{n}+b_{n}-b_{i+1}}\frac{a_{n}}{b_{n+1}}.
\end{eqnarray*}

We recall that the asymptotic behavior of the probability $\mathbf{P}\left(
\mathcal{A}_{i}(n)\right) $ as $n\rightarrow \infty $ and $i$ varies with $n$
in an appropriate way is described by the following theorem (see \cite%
{SV2021}).

\begin{theorem}
\label{T_total}If Hypotheses A1--A2 are valid then

1) for any fixed $N$%
\begin{equation*}
\lim_{n\rightarrow \infty }n^{1/2}\mathbf{P}\left( \mathcal{A}%
_{n-N}(n)\right) =r_{N}\in \left( 0,\infty \right) ;  \label{AsymEnd}
\end{equation*}

2) for any fixed $i$%
\begin{equation}
\lim_{n\rightarrow \infty }n^{3/2}\mathbf{P}\left( \mathcal{A}_{i}(n)\right)
=w_{i}\in \left( 0,\infty \right) ;  \label{AsumBegin}
\end{equation}

3) if, in addition, Hypothesis A3 is valid then
\begin{equation}
\lim_{\min \left( i,n-i\right) \rightarrow \infty }i^{1/2}\left( n-i\right)
^{3/2}\mathbf{P}\left( \mathcal{A}_{i}(n)\right) =K\in \left( 0,\infty
\right),  \label{AsymInterm}
\end{equation}
where an expression for $K=K(\infty)$ is given by formula (\ref%
{ExplicK}) below.
\end{theorem}

The next statement is a particular case of a theorem established in \cite{VD2020} (see \cite{GL2001} for a previous result, more general than what is needed in our case).

\begin{lemma}
\label{L_Guiv0} Let $h:$ $[0,\infty )\times \lbrack 0,\infty )\rightarrow
\lbrack 0,\infty )$ be a nonnegative continuous function not identically
equal zero and such that there exist constants $0<\lambda ~$\ and $C>0$ such
that
\begin{equation*}
\ h(x,y)\leq \frac{C}{\left( 1+x+y\right) ^{\lambda }}  \label{Estgh}
\end{equation*}%
for all $x\geq 0,y\geq 0$.

If the distribution of $X$ is nonlattice, $\mathbf{E}X=0$ and $varX\in
(0,\infty )$ then there exists a positive constant $K\left( h\right) $ such
that
\begin{equation*}
\lim_{n\rightarrow \infty }n^{1/2}\mathbf{E}\left[ h(a_{n},B_{1,n})\right]
=K(h).
\end{equation*}
\end{lemma}

With these results in hands, we will now be able to prove our main result,
Theorem \ref{T_conditional}.

\section{The case $i=n-N$}

\label{sec_nN}

Let $N$ be a fixed integer and consider the case $i=n-N$. Taking the
expectation with respect to the $\sigma $-algebra generated by the sequence $%
F_{i},F_{i+1},...,F_{n}$ and making the changes $F_{j}\rightarrow \tilde{F}%
_{j-i}$ for $j=i,\ldots ,n$ we write%
\begin{eqnarray*}
\mathbf{E}\left[ \mathcal{H}_{i,n}(s)\right] &=&\mathbf{E}\left[ \left(
1-F_{i,n}(s)\right)
\prod_{j=i+1}^{n-1}F_{j,n}(0)\prod_{j=0}^{i-1}F_{j,n-N}(F_{i,n}(0))\right] \\
&=&\mathbf{E}\left[ \left( 1-\tilde{F}_{0,N}(s)\right) \prod_{j=1}^{N}\tilde{%
F}_{0,j}(0)\prod_{j=0}^{i-1}F_{j,i}(\tilde{F}_{0,N}(0))\right] \\
&=&\mathbf{E}\left[ \frac{\left( 1-\tilde{F}_{0,N}(s)\right) }{%
\left( 1-\tilde{F}_{0,N}(0)\right) }\prod_{j=1}^{N}\tilde{F}_{0,j}(0)\Psi_{i}\left( \tilde{F}_{0,N}(0)\right) \right] ,
\end{eqnarray*}%
where%
\begin{equation*}
\Psi_{i}\left( z\right) :=(1-z)\mathbf{E}%
\left[ \prod_{j=0}^{i-1}F_{j,i}(z)\right] .
\end{equation*}

Observe that as the environments are i.i.d.,
\begin{equation*}
\prod_{k=1}^{i-1}F_{k,i}(z)\overset{d}{=}\prod_{k=1}^{i-1}F_{k,0}(z).
\end{equation*}%
Now from Lemma 2 in \cite{DLVZ19}
\begin{equation*}
\prod_{k=1}^{i-1}F_{k,0}(z)=\frac{(1-z)^{-1}}{(1-z)^{-1}+%
\sum_{k=1}^{i-1}e^{-S_{k}}}=\frac{(1-z)^{-1}}{%
(1-z)^{-1}+a_{i-1}+B_{1,i-1}}.
\end{equation*}%
By a direct application of Lemma \ref{L_Guiv0} with $\lambda =1$ we obtain that for any $z\in \lbrack 0,1)$,
\begin{equation*}
\lim_{i\rightarrow \infty }\sqrt{i}\Psi_{i}(z):=\psi (z)\in
(0,\infty ).
\end{equation*}%
Note also that in view of Theorem 1.1 in \cite{GK1999}
\begin{equation}
\sqrt{i}\Psi_{i}(z)\leq \sqrt{i}\mathbf{E}\left[ \frac{1}{%
1+\sum_{k=1}^{i-1}e^{-S_{k}}}\right] =\sqrt{i}\mathbf{P}\left(
Z_{i}>0\right) \leq C.  \label{Domin1}
\end{equation}

\textbf{Proof of point 1) of Theorem \ref{T_conditional}}. By conditioning
with respect to the piece of the environment $F_{1},...,F_{n}$, we may check that, for
all $s\in \lbrack 0,1]$%
\begin{equation}
\mathbf{E}\left[ s^{Z_{i,n}}|\mathcal{A}_{i}(n)\right] =1-\frac{\mathbf{E}%
\left[ 1-s^{Z_{i,n}};\mathcal{A}_{i}(n)\right] }{\mathbf{P}\left( \mathcal{A}%
_{i}(n)\right) }=1-\frac{\mathbf{E}\left[ \mathcal{H}_{i,n}(s)\right] }{%
\mathbf{P}\left( \mathcal{A}_{i}(n)\right) }.  \label{YaglomDiscrete}
\end{equation}%
From (\ref{Domin1}), we get
\begin{equation*}
\sqrt{i}\frac{(1-\tilde{F}_{0,N}(s))}{(1-\tilde{F}_{0,N}(0))}
\prod_{k=1}^{N-1}\tilde{F}_{k,N}(0)\Psi_{i}(\tilde{F}%
_{0,N}(0))\leq \sqrt{i}\Psi_{i}(\tilde{F}_{0,N}(0))\leq C.
\end{equation*}%
Applying the dominated convergence theorem we conclude that, for $i=n-N$%
\begin{eqnarray*}
\lim_{n\rightarrow \infty }\frac{\mathbf{E}\left[ \mathcal{H}_{i,n}(s)\right]
}{\mathbf{P}\left( \mathcal{A}_{i}(n)\right) } &=&\frac{\lim_{n\rightarrow
\infty }n^{1/2}\mathbf{E}\left[ \mathcal{H}_{i,n}(s)\right] }{%
\lim_{n\rightarrow \infty }n^{1/2}\mathbf{E}\left[ \mathcal{H}_{i,n}(0)%
\right] } \\
&=&{\frac{\lim_{i\rightarrow \infty }i^{1/2}\mathbf{E}\left[
\frac{1-\tilde{F}_{0,N}(s)}{{1-\tilde{F}_{0,N}(0))}}\prod_{k=1}^{N}\tilde{F%
}_{0,j}(0)\Psi_{i}\left( \tilde{F}_{0,N}(0)\right) \right] }{%
\lim_{i\rightarrow \infty }i^{1/2}\mathbf{E}\left[ \prod_{k=1}^{N}\tilde{F}%
_{0,j}(0)\Psi_{i}\left( \tilde{F}_{0,N}(0)\right) \right] }}
\\
&=&\frac{\mathbf{E}\left[ \frac{1-\tilde{F}_{0,N}(s)}{{1-\tilde{%
F}_{0,N}(0))}} \prod_{k=1}^{{N}}\tilde{F}_{k,N}(0)\psi (%
\tilde{F}_{0,N}(0))\right] }{\mathbf{E}\left[ \prod_{k=1}^{{N}}\tilde{F}%
_{k,N}(0)\psi (\tilde{F}_{0,N}(0))\right] }=:\phi
_{N}(s).
\end{eqnarray*}%
Let
\begin{equation*}
\Theta _{N}(s)=\mathbf{E}\left[ s^{\vartheta _{N}}\right] :=\lim_{n%
\rightarrow \infty }\mathbf{E}\left[ s^{Z_{i,n}}|\mathcal{A}_{i}(n)\right]
=1-\phi _{N}(s).
\end{equation*}%
Since $N$ is fixed, it follows from the previous relation that $\Theta
_{N}(0)=0,\, \Theta _{N}(1)=1$ and $\Theta _{N}(s)\in (0,1)$ for all $s\in
(0,1) $. Hence, $\vartheta _{N}$ is a proper nondegenerate random variable.

This proves (\ref{YagDiscrete}).

\section{The case $\,\min (i,n-i)\rightarrow \infty $}

\label{sec_ini_to_infty}

We now consider the case $\min (i,n-i)\rightarrow \infty $. Introduce the
running maximum and minimum of the associated random walk $\mathbf{S}:$
\begin{equation*}
M_{n}:=\max \left( S_{1},...,S_{n}\right) ,\quad L_{n}:=\min \left(
S_{0},S_{1},...,S_{n}\right)  \label{def_LM}
\end{equation*}%
and denote by
\begin{equation*}
\tau (n):=\min \{0\leq k\leq n:S_{k}=L_{n}\}  \label{def_tau(n)}
\end{equation*}%
the moment of the first random walk minimum up to time $n$.

It is known that if Hypothesis A2 is valid then (see, for instance, \cite%
{F08}, Ch.XII, Sec 7, Theorem 1a) there exist positive constants $K_{1}$ and
$K_{2}$ such that, as $n\rightarrow \infty $%
\begin{equation}
\mathbf{P}\left( \tau (n)=n\right) =\mathbf{P}(M_{n}<0)\sim
K_{1}n^{-1/2},\quad \mathbf{P}(L_{n}\geq 0)\sim K_{2}n^{-1/2},
\label{AsymPositive}
\end{equation}%
and (see, for instance, Proposition 2.1 in \cite{ABKV}) there
exist positive constants $K_{3}$ and $K_{4}$ such that, as $n\rightarrow
\infty $
\begin{equation}
\mathbf{E}\left[ e^{S_{n}};\tau (n)=n\right] \sim K_{3}n^{-3/2},\quad
\mathbf{E}\left[ e^{-S_{n}};L_{n}\geq 0\right] \sim K_{4}n^{-3/2}.
\label{AsymConditional}
\end{equation}

The proof of the third statement of Theorem \ref{T_conditional} is based on
two changes of measure performed by means of the right-continuous functions $%
U:\mathbf{R}$ $\rightarrow \lbrack 0,\infty )$ and $V:\mathbf{R}$ $%
\rightarrow \lbrack 0,\infty )$ specified by
\begin{equation*}
U(x):=I\{x\geq 0\}+\sum_{n=1}^{\infty }\mathbf{P}\left( S_{n}\geq -x,M_{n}<0\right),  \label{DefU}
\end{equation*}%
\begin{equation*}
V(x):=I\{x<0\}+\sum_{n=1}^{\infty }\mathbf{P}\left( S_{n}<-x,L_{n}\geq 0\right),
  \label{DefV}
\end{equation*}
where $I\{\mathcal{A}\}$ is the indicator of the set $\mathcal{A}.$

It is known (see, for instance, \cite{4h} and \cite{ABKV}) that for any oscillating random walk
\begin{equation}
\mathbf{E}\left[ U(x+X);X+x\geq 0\right] =U(x),\quad x\geq 0,  \label{Mes1}
\end{equation}%
and
\begin{equation}
\mathbf{E}\left[ V(x+X);X+x<0\right] =V(x),\quad x\leq 0.  \label{Mes2}
\end{equation}

Let $\mathcal{E}=\left\{ F_{n},n\in \mathbf{N}\right\} $ be a random
environment and let $\mathcal{F}_n$ be the $\sigma $%
-algebra of events generated by the random variables $F_{1},F_{2},...,F_{n}$
and the sequence $Y_{0},Y_{1},...,Y_{n}$. The sequence of these $\sigma $%
-algebras forms a filtration $\mathfrak{F}$. Using the martingale property (%
\ref{Mes1})-(\ref{Mes2}) of $U,V$ one can introduce (see, for instance,
\cite{KV2017}, Chapter 7) a sequence of probability measures $%
\{\mathbf{P}_{(n)}^{+},n\geq 1\}$ on the $\bigvee_{n\geq1}$
by means of the equalities
\begin{equation*}
d\mathbf{P}_{(n)}^{+}(\mathcal{A}):=U(S_{n})I\left\{ L_{n}\geq 0\right\} d\mathbf{P}(\mathcal{A}),\quad A\in\mathcal{F}_{n} .
\label{DefPplus}
\end{equation*}%%
Using this definition and Kolmogorov's extension theorem one can specify on
a suitable probability space a probability measure $\mathbf{P}^{+}$ on $\bigvee_{n\geq1}$ such that
\begin{equation}
\mathbf{P}^{+}|\mathcal{F}_{n}=\mathbf{P}_{(n)}^{+},\ n\in \mathbf{N}.
\label{DefMeasures}
\end{equation}

We write $\mathbf{P}_{x}$ and $\mathbf{E}_{x}$ for the corresponding
probability measures and expectations if $S_{0}=x$. Thus, $\mathbf{P}=%
\mathbf{P}_{0}.$ With this notation, (\ref{DefMeasures}) may be rewritten as
follows: for every $\mathcal{F}_{n}$-measurable random variable $O_{n}$ such that $\mathbf{E}_{x}\left[
O_{n}U(S_{n});L_{n}\geq 0\right]<\infty ,\ x\geq 0,$
\begin{equation*}
\mathbf{E}_{x}^{+}\left[ O_{n}\right] :=\frac{1}{U(x)}\mathbf{E}_{x}\left[
O_{n}U(S_{n});L_{n}\geq 0\right] ,\ x\geq 0.  \label{DefEplus}
\end{equation*}

Similarly, $V$ gives rise to probability measures $\mathbf{P}_{x}^{-},x\leq
0 $, which can be defined via:
\begin{equation*}
\mathbf{E}_{x}^{-}\left[ O_{n}\right] :=\frac{1}{V(x)}\mathbf{E}_{x}\left[
O_{n}V(S_{n});M_{n}<0\right] ,\ x\leq 0.  \label{DefEminus}
\end{equation*}

By means of the measures $\mathbf{P}_{x}^{+}$ and $\mathbf{P}_{x}^{-}$, we
investigate the limit behavior of certain conditional distributions.

For $\lambda >0$, let $\mu _{\lambda }$ and $\nu _{\lambda }$ be the
probability measures on $[0,+\infty )$ and $(-\infty ,0):$
\begin{equation*}
\mu _{\lambda }(dz)\ :=\ c_{1}e^{-\lambda z}U(z)1_{\{z\geq 0\}}\,dz\ ,\quad
\nu _{\lambda }(dz)\ :=\ c_{2}e^{\lambda z}V(z)1_{\{z<0\}}\,dz
\label{DefSubcrit}
\end{equation*}%
with
\begin{equation*}
c_{1}^{-1}=c_{1\lambda }^{-1}:=\int_{0}^{\infty }e^{-\lambda
z}U(z)\,dz,\quad c_{2}^{-1}=c_{2\lambda }^{-1}:=\int_{-\infty
}^{0}e^{\lambda z}V(z)\,dz.  \label{def_ci}
\end{equation*}

The next three lemmas are proven in \cite{SV2021} and are natural variations
of Lemmas 7.3 and 7.5 in \cite{KV2017}, Chapter 7. We recall
them for the sake of readability. We use the agreement $\delta n:=\lfloor
\delta n\rfloor $ for $0<\delta <1$ in their formulations and write below $%
\int^{(m)}$ for $m$ times integration $\int \cdots \int $.

Let $\mathcal{G},\mathcal{H},\mathcal{T}$ be three Euclidian (or Polish)
spaces. For each $r\in \mathbf{N}$ consider three functions $%
g_{r}:\Delta ^{r\delta }\rightarrow \mathcal{G},$ $h_{r}:\Delta ^{r\delta
}\rightarrow \mathcal{H}$, and $t_{r}:\Delta ^{r}\rightarrow \mathcal{T}$,
measurable with respect to the corresponding $\sigma-$algebras of Borel sets. \textbf{We assume in the three lemmas to follow that Hypotheses A2 and A3 hold.}

\begin{lemma}[Lemma 9 in \protect\cite{SV2021}]
\label{L_pr4New} Let $G_{r}:=g_{r}(F_{1},\ldots ,F_{\delta r})$ , $r\in
\mathbf{N}$ be random variables with values in $\mathcal{G}$ such that, as $%
r\rightarrow \infty $
\begin{equation*}
G_{r}\ \rightarrow \ G_{\infty }\quad \mathbf{P}^{+}-a.s.
\end{equation*}%
for some $\mathcal{G}$-valued random variable $G_{\infty }$. Also let $%
H_{r}:=h_{r}(F_{1},\ldots ,F_{\delta r})$, $r\in \mathbf{N}$, be random
variables with values in $\mathcal{H}$ such that, as $r\rightarrow \infty $
\begin{equation*}
H_{r}\ \rightarrow \ H_{\infty }\quad \mathbf{P}_{x}^{-}-a.s.
\end{equation*}%
for all $x\leq 0$ and some $\mathcal{H}$-valued random variable $H_{\infty }$%
. Denote
\begin{equation*}
\tilde{H}_{r}:=h_{r}(F_{r},...F_{r-\delta r+1}).
\end{equation*}%
Let $T_{r}:=t_{r}(F_{1},\ldots ,F_{r}),$ $r\in \mathbf{N}$ be random
variables with values in $\mathcal{T}$ such that, as $r\rightarrow \infty $
\begin{equation*}
T_{r}\ \rightarrow \ T_{\infty }\quad \mathbf{P}_{x}^{+}-a.s.
\end{equation*}%
for all $x\geq 0$ and some $\mathcal{T}$-valued random variable $T_{\infty }$%
. Denote $\hat{T}_{n-r}:=t_{n-r}(F_{r+1},\ldots ,F_{n})$ for $r\leq n$. Then
for $\lambda >0$ and any bounded continuous function $\varphi :\mathcal{G}%
\times \mathcal{H}\times \mathbf{R\times }\mathcal{T}\rightarrow \mathbf{R}$%
, as $\min \left( r,n-r\right) \rightarrow \infty ,$
\begin{eqnarray*}
&&\frac{\mathbf{E}[\varphi (G_{r},\tilde{H}_{r},S_{r};\hat{T}%
_{n-r})e^{-\lambda S_{r}}\;;\;L_{n}\geq 0]}{\mathbf{E}[e^{-\lambda
S_{r}};L_{r}\geq 0]\ \mathbf{P}\left( L_{n-r}\geq 0\right) } \\
&&\rightarrow\int^{(4)}U(-z)\varphi (u,v,-z,t)\mathbf{P}^{+}\left(
G_{\infty }\in du\right) \mathbf{P}_{z}^{-}\left( H_{\infty }\in dv\right)
\mathbf{P}_{-z}^{+}\left( T_{\infty }\in dt\right) \nu _{\lambda }(dz)\ .
\end{eqnarray*}
\end{lemma}

\begin{lemma}[Lemma 10 in \protect\cite{SV2021}]
\label{L_pr3New}Let $G_{r},H_{r},\tilde{H}_{r},T_{r},\hat{T}_{r}$, $r\in
\mathbf{N}$, be as in Lemma~\ref{L_pr4New}, now fulfilling, as $r\rightarrow
\infty $
\begin{equation*}
G_{r}\ \rightarrow \ G_{\infty }\quad \mathbf{P}_{x}^{+}-a.s.\,\forall x\geq
0,\quad H_{r}\rightarrow H_{\infty }\quad \mathbf{P}^{-}-a.s.,\quad
T_{r}\rightarrow T_{\infty }\quad \mathbf{P}^{+}-a.s.
\end{equation*}%
Then, for $\lambda >0$ and any bounded continuous function $\varphi :%
\mathcal{G}\times \mathcal{H}\times \mathbf{R\times }\mathcal{T}\rightarrow
\mathbf{R}$, as $\min \left( r,n-r\right) \rightarrow \infty $
\begin{eqnarray*}
&&\frac{\mathbf{E}[\varphi (G_{r},\tilde{H}_{r},S_{r};\hat{T}%
_{n-r})e^{\lambda S_{r}}\;;\;\tau (n)=r]}{\mathbf{E}[e^{\lambda S_{r}};\tau
(r)=r]\mathbf{P}\left( L_{n-r}\geq 0\right) } \\
&\rightarrow &\int^{(4)}\varphi (u,v,-z,t)\,\mathbf{P}_{z}^{+}\left(
G_{\infty }\in du\right) \mathbf{P}^{-}\left( H_{\infty }\in dv\right)
\mathbf{P}^{+}\left( T_{\infty }\in dt\right) \mu _{\lambda }(dz)\ .
\end{eqnarray*}
\end{lemma}

\begin{lemma}[Lemma 11 in \protect\cite{SV2021}]
\bigskip \label{L_pr6New}Let $G_{r},H_{r},\tilde{H}_{r},T_{r},\hat{T}_{r}$, $%
r\in \mathbf{N}$, be as in Lemma~\ref{L_pr4New}, and
$$
H_{N,r}=h_{N,r}(F_{N},F_{N+1},...,F_{r\delta -1}),\, \tilde{H}%
_{r,N}=h_{N,r}(F_{r-N},F_{r-N-1},...,F_{r-\delta r+1})
$$
now fulfilling as $r\rightarrow \infty $
\begin{equation*}
G_{r}\ \rightarrow \ G_{\infty }\quad \mathbf{P}_{x}^{+}-a.s.,\,\forall
x\geq 0,\quad \left( H_{r},H_{N,r}\right) \rightarrow (H_{\infty
},H_{N,\infty })\quad \mathbf{P}^{-}-a.s.
\end{equation*}%
and $T_{r}\rightarrow T_{\infty }$ $\mathbf{P}^{+}$-a.s. Then, for $\lambda
>0$ and for any bounded continuous function $\varphi :\mathcal{G}\times
\mathcal{H}\times \mathcal{H}\times \mathbf{R\times }\mathcal{T}\rightarrow
\mathbf{R}$, as $\min \left( r,n-r\right) \rightarrow \infty $
\begin{eqnarray*}
&&\frac{\mathbf{E}[\varphi (G_{r},\tilde{H}_{r},\tilde{H}_{N,r},S_{r};\hat{T}%
_{n-r})e^{\lambda S_{r}}\;;\;\tau (n)=r]}{\mathbf{E}[e^{\lambda S_{r}};\tau
(r)=r]\mathbf{P}\left( L_{n-r}\geq 0\right) } \\
&\rightarrow &\int^{(4)}\varphi (u,v_{1},v_{2},-z,t)\,\mathbf{P}%
_{z}^{+}\left( G_{\infty }\in du\right) \\
&&\times \mathbf{P}^{-}\left( (H_{\infty },H_{\infty ,N})\in
(dv_{1},dv_{2})\right) \mathbf{P}^{+}\left( T_{\infty }\in dt\right) \mu
_{\lambda }(dz)\ .
\end{eqnarray*}
\end{lemma}

We will need one more statement related to the driftless random walks $%
\left\{ S_n, n\in \mathbf{N}\right\} $.

For a fixed positive integer $N\leq \min (j/2,n-j)$ set
\begin{equation}  \label{defK1K2}
\mathcal{K}_{1}=\mathcal{K}_{1}(j,N):=[N,j-N]\cap \mathbf{N},\ \mathcal{K}%
_{2}=\mathcal{K}_{2}(j,n,N):=[j+N,n]\cap \mathbf{N}.
\end{equation}

\begin{lemma}
\label{L_Negl}[see Lemmas 13 and 14, and Corollary 15 in \cite{SV2021}] If
\begin{equation*}
\mathbf{E}X_{i}=0,\quad \sigma ^{2}=VarX_{i}\in \left( 0,\infty \right) ,
\label{CondRW}
\end{equation*}%
then there exists a constant $C>0$ such that for all $n\geq j\geq 1$
\begin{equation}
\mathbf{E}\left[ e^{-S_{j}}e^{S_{\tau (j-1)}}e^{S_{\tau (n)}}\right] \leq
\frac{C}{j^{3/2}\left( n-j+1\right) ^{1/2}}.  \label{AboveRW}
\end{equation}%
and for every $\varepsilon >0$ there exists $N=N(\varepsilon )$ such that
\begin{equation*}
\mathbf{E}\left[ e^{-S_{j}}e^{S_{\tau (j-1)}}e^{S_{\tau (n)}};\tau (n)\in
\mathcal{K}_{1}\cup \mathcal{K}_{2}\right] \leq \frac{\varepsilon }{j^{3/2}%
\sqrt{n-j+1}}
\end{equation*}%
for all $j\geq j_{0}=j_{0}\left( \varepsilon \right) ,n\geq
n_{0}=n_{0}\left( \varepsilon \right) $.
\end{lemma}

\bigskip We now prove point 3) of Theorem \ref{T_conditional}. Recall the
definition of $a_{i,n}$ in \eqref{defainbin}, and consider the rescaled
process $Y_{i,n}=a_{i,n}Z_{i,n}$.
Let $\beta >0$. If we take $s=\exp \left\{ -\beta a_{i,n}\right\} $ in (\ref%
{YaglomDiscrete}), we get
\begin{equation}
\mathbf{E}\left[ e^{-\mathcal{\beta }Y_{i,n}}|\mathcal{A}_{i}(n)\right] =1-%
\frac{\mathbf{E}\left[ \mathcal{H}_{i,n}(\exp \left\{ -\mathcal{\beta }%
a_{i,n}\right\} )\right] }{\mathbf{P}\left( \mathcal{A}_{i}(n)\right) }.
\label{YaglomForm}
\end{equation}

Note that in view of Lemma \ref{L_fractional}
\begin{equation*}
\mathcal{H}_{i,n}(\exp \left\{ -\mathcal{\beta }a_{i,n}\right\} )=\frac{1}{%
a_{i,n}/(1-\exp \left\{ -\mathcal{\beta }a_{i,n}\right\} )+b_{i,n}}\frac{%
a_{i,n}+b_{i,n}}{a_{i,n}+b_{i,n}-1}\frac{a_{n}}{a_n+b_n}.
\end{equation*}

Replacing $a_{i,n}$, $b_{i,n}$, $a_{n}$ and $b_{n}$ by their definition we
obtain,
\begin{eqnarray*}
&&\frac{1}{a_{i,n}/(1-\exp \left\{ -\mathcal{\beta }a_{i,n}\right\} )+b_{i,n}%
}\frac{a_{i,n}+b_{i,n}}{a_{i,n}+b_{i,n}-1}\frac{a_{n}}{a_{n}+b_{n}} \\
&=&\frac{1}{e^{S_{i}-S_{n}}/(1-\exp \left\{ -\mathcal{\beta }%
e^{S_{i}-S_{n}}\right\} )+\sum_{j=i}^{n-1}e^{S_{i}-S_{j}}} \\
&&\times \frac{e^{S_{i}-S_{n}}+\sum_{j=i}^{n-1}e^{S_{i}-S_{j}}}{%
e^{S_{i}-S_{n}}+\sum_{j=i+1}^{n-1}e^{S_{i}-S_{j}}}\frac{e^{-S_{n}}}{%
\sum_{j=0}^{n}e^{-S_{j}}} \\
&=&\frac{e^{S_{n}-S_{i}}}{1/(1-\exp \left\{ -\mathcal{\beta }%
e^{S_{i}-S_{n}}\right\} )+\sum_{j=i}^{n-1}e^{S_{n}-S_{j}}} \\
&&\times \frac{\sum_{j=i}^{n}e^{S_{n}-S_{j}}}{\sum_{j=i+1}^{n}e^{S_{n}-S_{j}}%
}\frac{1}{\sum_{j=0}^{n}e^{S_{n}-S_{j}}}.
\end{eqnarray*}%
Using the duality property of random walks (\cite{F08}, Ch.XII,
Section 2) we obtain the representation
\begin{multline}
\mathbf{E}\left[ \mathcal{H}_{i,n}(\exp \left\{ -\mathcal{\beta }%
a_{i,n}\right\} )\right] \\
=\mathbf{E}\left[ \frac{e^{S_{n-i}}}{1/(1-\exp \left\{ -\mathcal{\beta }%
e^{-S_{n-i}}\right\} )+\sum_{k=1}^{n-i}e^{S_{k}}}\frac{%
\sum_{k=0}^{n-i}e^{S_{k}}}{\sum_{k=0}^{n-i-1}e^{S_{k}}}\frac{1}{%
\sum_{k=0}^{n}e^{S_{k}}}\right] .  \label{NewInterpr}
\end{multline}

For the sake of readability, we now introduce the reflection
of the random walk $\mathbf{S}$,
\begin{equation*}
\bar{\mathbf{S}}:= \{\bar{S}_n , n \in \mathbf{N}_0\}= \{-S_n, n \in \mathbf{N}_0\}.
\end{equation*}
Functions and measures related to $\bar{\mathbf{S}}$ will be indicated with bars $%
\overset{\_}{\cdot }$. For instance, we write
\begin{equation*}
\bar{a}_{k}=e^{-\bar{S}_{k}},\ \bar{L}_{n}:=\min_{0\leq r\leq n}\bar{S}%
_{k},\ \bar{M}_{n}:=\max_{1\leq k\leq n}\bar{S}_{k},\ \bar{\tau}(n):=\min
\left\{ k\geq 0:\bar{S}_{k}=\bar{L}_{n}\right\}
\end{equation*}%
and
\begin{equation*}
\bar{U}(x):=I\{x\geq 0\}+\sum_{n=1}^{\infty }\mathbf{P}\left( \bar{S}_{n}\geq -x,\bar{M}%
_{n}<0\right)
\end{equation*}%
and specify the measure $\mathbf{\bar{P}}_{x}^{+}$ by the relation
\begin{equation*}
\mathbf{\bar{E}}_{x}^{+}\left[ O_{n}\right] :=\frac{1}{\bar{U}(x)}\mathbf{E}%
_{x}\left[ O_{n}\bar{U}(\bar{S}_{n});\bar{L}_{n}\geq 0\right] ,\ x\geq 0.
\end{equation*}

We also write $j$ instead of $n-i$ in the remaining part of the proofs. The
agreements above allow us to rewrite (\ref{NewInterpr}) as%
\begin{equation}
\mathbf{E}\left[ \mathcal{H}_{i,n}(\exp \left\{ -\mathcal{\beta }%
a_{i,n}\right\} )\right] =\mathbf{E}\left[ \mathcal{V}_{j,n}(\mathcal{\beta }%
)\right] ,  \label{intro_Nu}
\end{equation}%
where
\begin{equation}
\mathcal{V}_{j,n}(\mathcal{\beta }):=\frac{\bar{a}_{j}}{(1-\exp \left\{ -%
\mathcal{\beta }/\bar{a}_{j}\right\} )^{-1}+\bar{B}_{1,j+1}}\frac{\bar{b}%
_{j+1}}{\bar{b}_{j}}\frac{1}{\bar{b}_{n+1}},  \label{DefNu}
\end{equation}

It will be also convenient to consider
\begin{equation*}
\mathcal{V}_{j,n}(\mathcal{\infty }):=\frac{\bar{a}_{j}}{\bar{b}_{j}}\frac{1%
}{\bar{b}_{n+1}}.
\end{equation*}%
According to this agreement
\begin{equation*}
\mathbf{E}\left[ \mathcal{V}_{j,n}(\mathcal{\infty })\right] =\mathbf{E}%
\left[ \mathcal{H}_{i,n}(0)\right] =\mathbf{P}\left( \mathcal{A}%
_{i}(n)\right) .
\end{equation*}

It will be clear from the arguments to follow that all our estimates and
limiting expressions are valid for all $\beta \in (0,\infty ]$, i.e., $\beta
=\infty $ is included.

We fix some positive integer $N$, and recall the definition of $\mathcal{K}%
_1 $ and $\mathcal{K}_2$ in \eqref{defK1K2}. We have the decomposition
\begin{eqnarray}
\mathbf{E}\left[ \mathcal{V}_{j,n}(\mathcal{\beta })\right] &=&\mathbf{E}%
\left[ \mathcal{V}_{j,n}(\mathcal{\beta });\bar{\tau}(n)\in \mathcal{K}%
_{1}\cup \mathcal{K}_{2}\right]  \notag \\
&&+\mathbf{E}\left[ \mathcal{V}_{j,n}(\mathcal{\beta });\bar{\tau}(n)<N%
\right] +\mathbf{E}\left[ \mathcal{V}_{j,n}(\mathcal{\beta });\bar{\tau}%
(n)\in (j-N,j]\right]  \notag \\
&&+\mathbf{E}\left[ \mathcal{V}_{j,n}(\mathcal{\beta });\bar{\tau}(n)\in
(j,j+N)\right] .  \label{DecompYagl}
\end{eqnarray}%
Since
\begin{equation*}
\mathbf{E}\left[ \mathcal{V}_{j,n}(\mathcal{\beta })\Bigg|\mathbf{\bar{S}}%
\right] \leq \mathbf{E}\left[ \mathcal{V}_{j,n}(\infty )\Bigg|\mathbf{\bar{S}%
}\right] \overset{d}{=}\mathbf{P}\left( \mathcal{A}_{i}(n)|\mathbf{S}\right)
\end{equation*}%
and
\begin{equation*}
\mathbf{E}\left[ \mathcal{V}_{j,n}(\mathcal{\beta })\right] \leq \mathbf{P}%
\left( \mathcal{A}_{i}(n)\right) ,\quad \beta \in (0,\infty ],
\end{equation*}%
it follows from Lemma 14 in \cite{SV2021} that for all $n\geq j+1\geq 2$ and
for any $\varepsilon >0$ there exists $N\left( \varepsilon \right) $ such
that, for all $N\geq N(\varepsilon )$ and all sufficiently large $j$ and $%
n-j $
\begin{equation*}
\mathbf{E}\left[ \mathcal{V}_{j,n}(\mathcal{\beta });\bar{\tau}(n)\in
\mathcal{K}_{1}\cup \mathcal{K}_{2}\right] \leq \frac{\varepsilon }{%
j^{3/2}\left( n-j\right) ^{1/2}},\ \beta \in (0,\infty ].
\end{equation*}

We now focus on the last three summands in the right-hand side of (\ref{DecompYagl}).

We write
\begin{equation*}
\bar{b}_{n}=\bar{b}_{k}+\bar{a}_{k}\bar{b}_{k,n}
\end{equation*}
and use below the equality $\bar{b}_{k,n}\overset{d}{=}\bar{b}_{n-k}$ and
the independence of $\left( \bar{b}_{k},\bar{a}_{k}\right) $ and $\bar{b}%
_{k,n}$ many times.

1) We first evaluate $\mathbf{E}\left[ \mathcal{V}_{j,n}(\mathcal{\beta });%
\bar{\tau}(n)<N\right] $. Let $k<N<j$. By conditioning on the trajectory of
the random walk $\mathbf{\bar{S}}$ until time $k$ we obtain
\begin{equation*}  \label{Nutaunk}
\mathbf{E}\left[ \mathcal{V}_{j,n}(\mathcal{\beta });\bar{\tau}(n)=k\right] =%
\mathbf{E}\left[ e^{-\bar{S}_{k}}\Psi _{j-k,n-k}\left( \mathcal{\beta };e^{-%
\bar{S}_{k}},\bar{b}_{k}\right) ;\bar{\tau}(k)=k\right]
\end{equation*}
where
\begin{multline*}
\Psi _{j,n}\left( \mathcal{\beta };u,q\right) := \\
\mathbf{E}\left[ e^{-\bar{S}_{j}}\frac{1}{(1-\exp \left\{ -\mathcal{\beta }/u%
\bar{a}_{j}\right\} )^{-1}+q-1+u\bar{b}_{j+1}}\frac{q+u\bar{b}_{j+1}}{q+u%
\bar{b}_{j}}\frac{1}{q+u\bar{b}_{n+1}};\bar{L}_{n}\geq 0\right]
\label{ExpressLambda}
\end{multline*}
for $\beta \in (0,\infty )$ and
\begin{equation*}
\Psi _{j,n}\left( \mathcal{\infty };u,q\right) :=\mathbf{E}\left[ e^{-\bar{S}%
_{j}}\frac{1}{q+u\bar{b}_{j}}\frac{1}{q+u\bar{b}_{n+1}};\bar{L}_{n}\geq 0%
\right] .  \label{ExpressLambda_infty}
\end{equation*}

Set $t=\left[ j/2\right] $ and denote for $w\geq t+2$
\begin{eqnarray}
&&G_{t}:=\sum_{r=0}^{t}e^{-\bar{S}_{r}},\ \tilde{H}%
_{t+1,w-1}:=\sum_{r=t+1}^{w-1}e^{\bar{S}_{w}-%
\bar{S}_{r}}\overset{d}{=}\sum_{r=1}^{w-t-1}e^{\bar{S}%
_{r}}=:H_{1,w-t},  \notag \\
&& \hat{T}_{j,n} :=\sum_{r=j}^{n}e^{\bar{S}_{j}-\bar{S}_{r}}\overset{d}{=}%
\sum_{r=0}^{n-j}e^{-\bar{S}_{r}}=:T_{n-j}.  \label{NewT}
\end{eqnarray}%
Observe that%
\begin{equation*}
\bar{b}_{j+1}=G_{t}+e^{-\bar{S}_{j}}\tilde{H}%
_{t+1,j-1}+e^{-\bar{S}_{j}}=\bar{b}_{j}+e^{-\bar{S}_{j}}  \label{Reprb_j+1}
\end{equation*}%
and%
\begin{equation*}
\bar{b}_{n+1}=G_{t}+e^{-\bar{S}_{j}}\tilde{H}%
_{t+1,j-1}+e^{-\bar{S}_{j}}\ \hat{T}_{j,n}.  \label{Reprb_n+1}
\end{equation*}%
Using these equalities we may check that
\begin{equation*}
\Psi _{j,n}\left( \mathcal{\beta };u,q\right) =\mathbf{E}\left[ e^{-\bar{S}%
_{j}}\varphi _{\mathcal{\beta },u,q}(G_{t},\tilde{H}%
_{t+1,j-1},\bar{S}_{j};\hat{T}_{j,n});\bar{L}_{n}\geq 0\right] ,
\end{equation*}%
where
\begin{eqnarray*}
\varphi _{\mathcal{\beta },u,q}(g,h,s;t) &:&=\frac{1}{(1-\exp \left\{ -%
\mathcal{\beta }/ue^{-s}\right\} )^{-1}+q-1+u\left( g+e^{-s}(h+1\right) )} \\
&&\times \frac{q+u\left( g+e^{-s}(h+1)\right) }{q+u\left( g+e^{-s}h\right) }%
\frac{1}{q+u\left( g+e^{-s}(h+t\right) )}.
\end{eqnarray*}%
Clearly, $\varphi _{\mathcal{\beta },u,q}(g,h,s;t)$ is a continuous function
for $\beta\geq 0,\, q\geq 1$ and $u,g,h,s,t\geq 0$ and, by monotonicity in $%
\beta $%
\begin{eqnarray}
&&\varphi _{\mathcal{\beta },u,q}(g,h,s;t)\leq \varphi _{\mathcal{\infty }%
,u,q}(g,h,s;t)\notag\\
&&\qquad\qquad\qquad\quad:=\frac{1}{q+u\left( g+e^{-s}h\right) }\frac{1}{q+u\left(
g+e^{-s}(h+t\right) )}\leq 1.  \label{EstPhi_beta}
\end{eqnarray}

Since the random walk $\mathbf{\bar{S}}$ is driftless and $var\bar{X}%
_{i}=varX_{i}$, it follows from Lemma 5.5 in \cite{GV2017} that, as $\min (j,n-j)\rightarrow \infty $,
\begin{eqnarray*}
G_{t} &\rightarrow &G_{\infty }:=\sum_{r=0}^{\infty }e^{-\bar{%
S}_{r}}\quad \mathbf{\bar{P}}^{+}-a.s.,  \notag \\
H_{1,j-t} &\rightarrow &H_{1,\infty }:=\sum_{r=1}^{\infty }e^{\bar{S}%
_{r}}\quad \;\mathbf{\bar{P}}_{x}^{-}-a.s.,\quad \forall x\leq 0,  \notag \\
T_{n-j} &\rightarrow &T_{\infty }:=\sum_{r=0}^{\infty }e^{-\bar{S}_{r}}\quad
\mathbf{\bar{P}}_{x}^{+}-a.s.,\quad \forall x\geq 0.  \label{ShortT_i}
\end{eqnarray*}%
Estimate (\ref{EstPhi_beta}) allows us to apply Lemma~\ref{L_pr4New} and to
obtain that
\begin{equation}
\lim_{\min (j,n-j)\rightarrow \infty }\frac{\Psi _{j,n}\left( \mathcal{\beta
};u,q\right) }{\mathbf{E}[e^{-\bar{S}_{j}};\bar{L}_{j}\geq 0]\mathbf{P}%
\left( \bar{L}_{n-j}\geq 0\right) \ }=\Psi _{\infty }\left( \mathcal{\beta }%
;u,v\right)  \label{Psi_j,n}
\end{equation}%
exists for each fixed tuple $(\mathcal{\beta };u,q)$, where, for $\beta \in
(0,\infty ]$%
\begin{eqnarray*}
&&\Psi _{\infty }\left( \mathcal{\beta };u,q\right):=\int^{(4)}\bar{U}%
(-s)\varphi _{\beta ,u,q}(g,h,-s,t)\mathbf{\bar{P}}^{+}\left( G_{\infty }\in
dg\right)  \notag \\
&&\times \mathbf{\bar{P}}_{s}^{-}\left( H_{1,\infty }\in dh\right) \mathbf{%
\bar{P}}^{+}_{-s}\left( T_{\infty }\in dt\right) \bar{\nu}%
_{1}(ds).  \label{Def_Psi_beta}
\end{eqnarray*}%
Since
\begin{equation*}
\lim_{\beta \downarrow 0}\varphi _{\mathcal{\beta },u,q}(g,h,s;t)=0
\end{equation*}%
and the inequality (\ref{EstPhi_beta}) is valid, we conclude by the
dominated convergence theorem that%
\begin{equation*}
\lim_{\beta \downarrow 0}\Psi _{\infty }\left( \mathcal{\beta };u,q\right)
=0.  \label{Psi_zero}
\end{equation*}%
Furthermore, setting
\begin{equation*}
C_{k}(\mathcal{\beta }):=\mathbf{E}\left[ e^{-\bar{S}_{k}}\Psi _{\infty
}\left( \mathcal{\beta };e^{-\bar{S}_{k}},\bar{B}_{1,k}\right) ;\bar{\tau}%
(k)=k\right]  \label{DefC_k}
\end{equation*}%
we have, again by monotonicity of $\Psi _{\infty }\left( \mathcal{\beta }%
;u,q\right) $ in $\beta $ and the dominated convergence theorem that
\begin{equation}
\lim_{\beta \downarrow 0}C_{k}(\mathcal{\beta })=0  \label{C_k_vanish}
\end{equation}%
and
\begin{equation}
\lim_{\beta \uparrow \infty }C_{k}(\mathcal{\beta })=\mathbf{E}\left[ e^{-%
\bar{S}_{k}}\Psi _{\infty }\left( \mathcal{\infty };e^{-\bar{S}_{k}},\bar{B}%
_{1,k}\right) ;\bar{\tau}(k)=k\right] =C_{k}(\infty )<\infty .
\label{C_k_infinite}
\end{equation}

Invoking the dominated convergence theorem once more we conclude by (\ref%
{Psi_j,n}) (\ref{AsymPositive}) and (\ref{AsymConditional}) that for any $%
k<N $, as $\min \left( j,n-j\right) \rightarrow \infty $%
\begin{eqnarray}
&&\mathbf{E}\left[ e^{-\bar{S}_{k}}\Psi _{j-k,n-k}\left( \mathcal{\beta }%
;e^{-\bar{S}_{k}},\bar{B}_{1,k}\right) ;\bar{\tau}(k)=k\right]  \notag \\
&\sim &\mathbf{E}\left[ e^{-\bar{S}_{k}}\Psi _{\infty }\left( \mathcal{\beta
};e^{-\bar{S}_{k}},\bar{B}_{1,k}\right) ;\bar{\tau}(k)=k\right] \mathbf{E}%
[e^{-\bar{S}_{j}};\bar{L}_{j}\geq 0]\mathbf{P}\left( \bar{L}_{n-j}\geq
0\right)  \notag \\
&&\qquad \qquad \qquad \qquad \sim \frac{K_{2}K_{4}C_{k}(\mathcal{\beta })}{%
j^{3/2}\left( n-j\right) ^{1/2}}.  \label{TermLeft}
\end{eqnarray}%
Recalling (\ref{intro_Nu}) and summing (\ref{TermLeft}) over the $k$'s in $%
[0,N-1]$, gives, as $\min (j,n-j)\rightarrow \infty $
\begin{equation}
j^{3/2}\left( n-j\right) ^{1/2}\mathbf{E}\left[ \mathcal{V}_{j,n}(\mathcal{%
\beta });\bar{\tau}(n)<N\right] \sim K_{2}K_{4}C(\mathcal{\beta },N)
\label{LeftPart1}
\end{equation}%
where%
\begin{equation}
C(\mathcal{\beta },N):=\sum_{k=0}^{N}\mathbf{E}\left[ e^{-\bar{S}_{k}}\Psi
_{\infty }\left( \mathcal{\beta };e^{-\bar{S}_{k}},\bar{B}_{1,k}\right) ;%
\bar{\tau}(k)=k\right] .  \label{LeftPart2}
\end{equation}

2) We now evaluate $\mathbf{E}\left[ \mathcal{V}_{j,n}(\mathcal{\beta });%
\bar{\tau}(n)\in (j-N,j]\right] $. To this aim we fix $1\leq k<N$, recall
that $t=\left[ j/2\right] $, as well as the definitions of $\tilde{H}$ and $%
\hat{T}$ in (\ref{NewT}). We also introduce
\begin{equation*}
\tilde{D}_{j,k}:=\bar{S}_{j-k}-\bar{S}_{j}\overset{d}{=}-\bar{S}_{k}.
\end{equation*}%
The same as before, as $\min (j,n-j)\rightarrow \infty $
\begin{equation*}
G_{t}\rightarrow G_{\infty }\quad \mathbf{\bar{P}}%
_{x}^{+}-a.s.\,\forall x\geq 0,\,H_{1,j-t-k-1}\rightarrow H_{1,\infty }\quad
\mathbf{\bar{P}}^{-}-a.s.,\,T_{0,n-j+k}\rightarrow T_{\infty }\quad \mathbf{%
\bar{P}}^{+}-a.s.
\end{equation*}

Besides, using \eqref{DefNu} and multiplying both numerator and denominator
by $e^{2\bar{S}_{j-k}}$, we obtain
\begin{eqnarray*}
&&\mathbf{E}\left[ \mathcal{V}_{j,n}(\mathcal{\beta });\bar{\tau}(n)=j-k%
\right] \\
&&=\mathbf{E}\left[ e^{\bar{S}_{j-k}}\frac{e^{\bar{S}_{j-k}-\bar{S}_{j}}}{e^{%
\bar{S}_{j-k}}\left( 1-\exp \left\{ -\mathcal{\beta }e^{\bar{S}_{j}}\right\}
\right) ^{-1}+e^{\bar{S}_{j-k}}\bar{B}_{1,j+1}}\frac{\bar{b}_{j+1}}{\bar{b}%
_{j}}\frac{1}{e^{\bar{S}_{j-k}}\bar{b}_{n+1}};\bar{\tau}(n)=j-k\right]
\end{eqnarray*}%
Observing that%
\begin{eqnarray*}
e^{\bar{S}_{j-k}}\bar{B}_{1,j+1} &=&e^{\bar{S}_{j-k}}\sum_{r=1}^{t}e^{-\bar{S%
}_{r}}+\sum_{r=t+1}^{j-k-1}e^{\bar{S}_{j-k}-\bar{S}_{r}}+\sum_{r=j-k}^{j}e^{%
\bar{S}_{j-k}-\bar{S}_{r}} \\
&=&e^{\bar{S}_{j-k}}(G_{t}-1)+\tilde{H}_{t+1,j-k}+%
\hat{T}_{j-k,j}
\end{eqnarray*}%
we write
\begin{eqnarray*}
&&\frac{e^{\bar{S}_{j-k}-\bar{S}_{j}}}{e^{\bar{S}_{j-k}}\left( 1-\exp
\left\{ -\mathcal{\beta }e^{\bar{S}_{j}}\right\} \right) ^{-1}+e^{\bar{S}%
_{j-k}}\bar{B}_{1,j+1}} \\
&=&\frac{e^{\bar{S}_{j-k}-\bar{S}_{j}}}{e^{\bar{S}_{j-k}}\left( 1-\exp
\left\{ -\mathcal{\beta }e^{\bar{S}_{j}-\bar{S}_{j-k}}e^{\bar{S}%
_{j-k}}\right\} \right) ^{-1}+e^{\bar{S}_{j-k}}(G_{t}-1)+\tilde{H}_{t+1,j-k}+\hat{T}_{j-k,j}}
\end{eqnarray*}%
and, similarly,%
\begin{equation*}
e^{\bar{S}_{j-k}}\bar{b}_{j+1}=e^{\bar{S}_{j-k}}G_{t}+\tilde{H}_{t+1,j-k}+\hat{T}_{j-k,j},
\end{equation*}%
\begin{equation*}
e^{\bar{S}_{j-k}}\bar{b}_{j}=e^{\bar{S}_{j-k}}G_{t}+
\tilde{H}_{t+1,j-k}+\hat{T}_{j-k,j-1},
\end{equation*}%
\begin{equation*}
e^{\bar{S}_{j-k}}\bar{b}_{n+1}=e^{\bar{S}_{j-k}}G_{t}+\tilde{H}_{t+1,j-k}+\hat{T}_{j-k,n}.
\end{equation*}%
Thus,
\begin{eqnarray*}
&&\mathbf{E}\left[ \mathcal{V}_{j,n}(\mathcal{\beta });\bar{\tau}(n)=j-k%
\right] \\
&=&\mathbf{E}\left[ e^{\bar{S}_{j-k}}\varphi _{\beta }^{-}(G_{t%
},\tilde{H}_{t+1,j-k},\bar{S}_{j-k},\tilde{D}_{j,k},\hat{T}%
_{j-k,j},\hat{T}_{j-k,j-1},\hat{T}_{j-k,n});\bar{\tau}(n)=j-k\right] ,
\end{eqnarray*}%
where, for $\beta \in \left( 0,\infty \right) $
\begin{eqnarray*}
\varphi _{\beta }^{-}(g,h,s_{1},s_{2},t_{1},t_{2},t_{3}) &:&=\frac{e^{s_{2}}%
}{e^{s_{1}}(1-\exp \left\{ -\mathcal{\beta }e^{-s_{2}}e^{s_{1}}\right\}
)^{-1}+e^{s_{1}}(g-1)+h+t_{1}} \\
&&\times \frac{e^{s_{1}}g+h+t_{1}}{e^{s_{1}}g+h+t_{2}}\frac{1}{%
e^{s_{1}}g+h+t_{3}}
\end{eqnarray*}%
and, by monotonicity of the function with respect to $\beta$
\begin{eqnarray*}
\varphi _{\beta }^{-}(g,h,s_{1},s_{2},t_{1},t_{2},t_{3}) &\leq &\varphi
_{\infty }^{-}(g,h,s_{1},s_{2},t_{1},t_{2},t_{3})  \notag \\
&=&\frac{e^{s_{2}}}{e^{s_{1}}g+h+t_{2}}\frac{1}{e^{s_{1}}g+h+t_{3}}\leq 1
\label{Phi_minus_bound}
\end{eqnarray*}%
in the domain%
\begin{equation*}
\left\{ g\geq 1,h\geq 0,s_{1}\leq 0,s_{2}\leq 0,t_{1}\geq 0,t_{2}\geq
0,t_{3}\geq 0\right\} \cap \left\{ e^{s_{2}}\leq
e^{s_{1}}g+h+t_{2},e^{s_{1}}g+h+t_{3}\geq 1\right\} .
\end{equation*}%
The same as before, as $\min (j,n-j)\rightarrow \infty $
\begin{eqnarray*}
G_{t} &\rightarrow &G_{\infty }\quad \mathbf{\bar{P}}%
_{x}^{+}-a.s.\,\forall x\geq 0,\,H_{1,j}\rightarrow H_{1,\infty }\quad
\mathbf{\bar{P}}^{-}-a.s.,\, \\
(T_{k-1},T_{k},T_{n-j+k}) &\rightarrow &(T_{k-1},T_{k},T_{\infty })\quad
\mathbf{\bar{P}}^{+}-a.s.
\end{eqnarray*}%
Recalling Lemma \ref{L_pr6New} we see that, for each $k$ there exists a
constant $J_{-k}(\mathcal{\beta })\geq 0$ such that, as $\min
(j,n-j)\rightarrow \infty $
\begin{eqnarray}
\mathbf{E}\left[ \mathcal{V}_{j,n}(\mathcal{\beta });\bar{\tau}(n)=j-k\right]
&\sim &J_{-k}(\mathcal{\beta })\mathbf{E}[e^{\bar{S}_{j-k}};\bar{\tau}%
(j-k)=j-k]\mathbf{P}\left( \bar{L}_{n-j+k}\geq 0\right)  \notag \\
&\sim &\frac{K_{2}K_{3}J_{-k}(\mathcal{\beta })}{j^{3/2}\left( n-j\right)
^{1/2}},  \label{Term2a}
\end{eqnarray}%
where
\begin{eqnarray*}
&&J_{-k}(\mathcal{\beta }) :=\int^{(7)}\varphi _{\beta
}^{-}(g,h,-s_{1},s_{2},t_{1},t_{2},t_{3})\mathbf{\bar{P}}^{+}\left(
G_{\infty }\in dg\right) \mathbf{\bar{P}}_{s_{1}}^{-}\left( (H_{1,\infty },-%
\bar{S}_{k})\in d(h,s_{2})\right)  \notag \\
&&\qquad\qquad\qquad\quad\times \mathbf{\bar{P}}^{+}\left( (T_{k-1},T_{k},T_{\infty })\in
d(t_{1},t_{2},t_{3})\right) \bar{\mu}_{1}(ds_{1}).  \label{Jminus}
\end{eqnarray*}

Monotonicity of $\varphi _{\beta }^{-}(g,h,-s_{1},-s_{2},t_{1},t_{2},t_{3})$
in $\beta $ and the dominated convergence theorem show that, for each $k$%
\begin{equation}
\lim_{\beta \downarrow 0}J_{-k}(\mathcal{\beta })=0  \label{J_minus_vanish}
\end{equation}%
and%
\begin{equation}
\lim_{\beta \uparrow \infty }J_{-k}(\mathcal{\beta })=J_{-k}(\mathcal{\infty
})<\infty .  \label{J_minus_infinity}
\end{equation}%
Summing (\ref{Term2a}) over the $k$'s in $[1,N-1]$, gives, as $\min
(j,n-j)\rightarrow \infty $
\begin{multline}
j^{3/2}\left( n-j\right) ^{1/2}\mathbf{E}\left[ \mathcal{V}_{j,n}(\mathcal{%
\beta });\bar{\tau}(n)\in (j-N,j)\right] \\ \sim K_{2}K_{3} J_{-}(\mathcal{%
\beta },N)= K_{2}K_{3} \sum_{k=1}^{N-1
}J_{-k}(\mathcal{\beta }).  \label{LeftTerm2}
\end{multline}

3) We finally evaluate $\mathbf{E}\left[ \mathcal{V}_{j,n}(\mathcal{\beta });%
\bar{\tau}(n)\in (j,j+N)\right] $. As before, we fix $1\leq k\leq N$, set $t=%
\left[ j/2\right] $ and introduce a more complicated notation
\begin{eqnarray*}
&&\tilde{H}_{t,w,v} :=\sum_{r=t+1}^{w}e^{\bar{S%
}_{w+v}-\bar{S}_{r}}\overset{d}{=}\sum_{r=v}^{ w+v-t-1}e^{\bar{S}_{r}}=:H_{v,w+v-t}, \\
&& \tilde{C}_{j,k}:=\bar{S}_{j+k}-\bar{S}_{j}\overset{d}{=}%
\bar{S}_{k}.
\end{eqnarray*}

With this notation in view we write
\begin{eqnarray*}
&&\frac{e^{\bar{S}_{j+k}-\bar{S}_{j}}}{e^{\bar{S}_{j+k}}\left( 1-\exp
\left\{ -\mathcal{\beta }e^{\bar{S}_{j}}\right\} \right) ^{-1}+e^{\bar{S}%
_{j+k}}\bar{B}_{1,j+1}} \\
&=&\frac{e^{\bar{S}_{j+k}-\bar{S}_{j}}}{e^{\bar{S}_{j+k}}\left( 1-\exp
\left\{ -\mathcal{\beta }e^{-\tilde{C}_{j,k}}e^{\bar{S}%
_{j+k}}\right\} \right) ^{-1}+e^{\bar{S}_{j+k}}(G_{t}-1)+
\tilde{H}_{t,j,k}},
\end{eqnarray*}%
\begin{equation*}
e^{\bar{S}_{j+k}}\bar{b}_{j+1}=e^{\bar{S}_{j+k}}G_{t}+\tilde{H}_{t,j,k}=e^{\bar{S}_{j+k}}\bar{b}_{j}+e^{\bar{S}_{j+k}-%
\bar{S}_{j}}
\end{equation*}%
and%
\begin{equation*}
e^{\bar{S}_{j+k}}\bar{b}_{n+1}=e^{\bar{S}_{j+k}}G_{t}+\tilde{H}_{t,j+k-1,1}+\hat{T}_{j+k,n}.
\end{equation*}%
As a result we have, using again \eqref{DefNu} and multiplying both
numerator and denominator by $e^{2\bar{S}_{j+k}}$
\begin{eqnarray*}
&&\mathbf{E}\left[ \mathcal{V}_{j,n}(\mathcal{\beta });\bar{\tau}(n)=j+k%
\right] \\
&=&\mathbf{E}\left[ e^{\bar{S}_{j+k}}\frac{e^{\bar{S}_{j+k}-\bar{S}_{j}}}{e^{%
\bar{S}_{j+k}}\left( 1-\exp \left\{ -\mathcal{\beta }e^{\bar{S}_{j}}\right\}
\right) ^{-1}+e^{\bar{S}_{j+k}}\bar{B}_{1,j+1}}\frac{\bar{b}_{j+1}}{\bar{b}%
_{j}}\frac{1}{e^{\bar{S}_{j+k}}\bar{b}_{n+1}};\bar{\tau}(n)=j+k\right] \\
&=&\mathbf{E}\left[ e^{\bar{S}_{j+k}}\varphi _{\beta }(G_{t},
\tilde{H}_{t,j,k},\tilde{H}_{t,j-1,k+1},\tilde{H}_{t,j+k-1,1},
\bar{S}_{j+k},\tilde{C}_{j,k},\hat{T}_{j+k,n});\bar{\tau}(n)=j+k%
\right] ,
\end{eqnarray*}%
where
\begin{eqnarray*}
\varphi _{\beta }^{+}(g,h_{1},h_{2},h_{3},s_{1},s_{2},t) &:&=\frac{e^{s_{2}}%
}{e^{s_{1}}(1-\exp \left\{ -\mathcal{\beta }e^{-s_{2}}e^{s_{1}}\right\}
)^{-1}+e^{s_{1}}(g-1)+h_{1}} \\
&&\times \frac{e^{s_{1}}g+h_{1}}{e^{s_{1}}g+h_{2}}\frac{1}{e^{s_{1}}g+h_{3}+t%
}
\end{eqnarray*}%
and%
\begin{eqnarray*}
\varphi _{\beta }^{+}(g,h_{1},h_{2},h_{3},s_{1},s_{2},t) &\leq &\varphi
_{\infty }^{+}(g,h_{1},h_{2},h_{3},s_{1},s_{2},t)  \notag \\
&:&=\frac{e^{s_{2}}}{e^{s_{1}}g+h_{2}}\frac{1}{e^{s_{1}}g+h_{3}+t}\leq 1
\label{Phi_plus_bound}
\end{eqnarray*}%
in the domain%
\begin{equation*}
\left\{ g\geq 1,h_{2}\geq h_{1}\geq 0,h_{3}\geq 0,s_{1}\leq 0,s_{2}\leq
0,t\geq 0\right\} \cap \left\{ e^{s_{2}}\leq
e^{s_{1}}g+h_{2},e^{s_{1}}g+h_{3}+t\geq 1\right\} .
\end{equation*}

We know that for any fixed $k$, as $\min (j,n-j)\rightarrow \infty $%
\begin{eqnarray*}
&&G_{t}\rightarrow G_{\infty }\quad \mathbf{\bar{P}}%
_{x}^{+}-a.s.,\quad \forall x\geq 0, \\
&&\left( H_{k,j-t-1},H_{k+1,j-t},H_{1,j+k-t},\bar{S}_{k}\right) \rightarrow
(H_{k,\infty },H_{k+1,\infty },H_{1,\infty },\bar{S}_{k})\;\;\mathbf{\bar{P}}%
^{-}-a.s., \\
&&T_{n-j-k}\rightarrow T_{\infty }\quad \mathbf{\bar{P}}^{+}-a.s.
\end{eqnarray*}%
Hence, using Lemma \ref{L_pr6New} and relations (\ref{AsymPositive})
and (\ref{AsymConditional}) we get that, for each fixed $k$ as $\min (j,n-j)\rightarrow \infty $%
\begin{eqnarray}
\mathbf{E}\left[ \mathcal{V}_{j,n}(\mathcal{\beta });\bar{\tau}(n)=j+k\right]
&\sim &J_{+k}(\mathcal{\beta })\mathbf{E}[e^{\bar{S}_{j+k}};\bar{\tau}%
(j+k)=j+k]\mathbf{P}\left( \bar{L}_{n-j-k}\geq 0\right)  \notag \\
&\sim &\frac{K_{2}K_{3}J_{+k}(\mathcal{\beta })}{j^{3/2}\left( n-j\right)
^{1/2}}.  \label{Term3}
\end{eqnarray}%
where
\begin{eqnarray*}
&&J_{+k}(\mathcal{\beta }):=\int^{(7)}\varphi _{\beta
}^{+}(g,h_{1},h_{2},h_{3},-s_{1},-s_{2},t_{1})\mathbf{\bar{P}}^{+}\left(
G_{\infty }\in dg\right) \notag \\
&&\times \mathbf{\bar{P}}^{+}\left( T_{\infty }\in
dt_{1}\right)\mathbf{\bar{P}}_{s_{1}}^{-}\left( (H_{k,\infty },H_{k+1,\infty
},H_{1,\infty },\bar{S}_{k})\in d(h_{1},h_{2},h_{3},s_{2}\right) )\bar{\mu}%
_{1}(ds_1).  \label{Jplus}
\end{eqnarray*}%
Similarly to the previous case, for each $k$
\begin{equation}
\lim_{\beta \downarrow 0}J_{+k}(\mathcal{\beta })=0.  \label{Jplus_vanish}
\end{equation}%
and%
\begin{equation}
\lim_{\beta \uparrow \infty }J_{+k}(\mathcal{\beta })=J_{+k}(\mathcal{\infty
})<\infty .  \label{J_plus_infinity}
\end{equation}%
Summing (\ref{Term3}) over the $k$'s in $[0,N-1]$, give, as $\min
(j,n-j)\rightarrow \infty $
\begin{equation}
\mathbf{E}\left[ \mathcal{V}_{j,n}(\mathcal{\beta });\bar{\tau}(n)\in
\lbrack j,j+N)\right] \sim \frac{J_{+}(\mathcal{\beta },N)}{j^{3/2}\left(
n-j\right) ^{1/2}}  \label{RightTerm3}
\end{equation}%
for some constant $J_{+}(\mathcal{\beta },N)$.

Combining (\ref{LeftPart1}), (\ref{LeftPart2}) and (\ref{RightTerm3}) with (%
\ref{LeftTerm2}) shows that, for any $\beta \in (0,\infty ]$
\begin{equation*}
\lim_{\min (j,n-j)\rightarrow \infty }j^{3/2}\left( n-j\right) ^{1/2}\mathbf{%
E}\left[ \mathcal{V}_{j,n}(\mathcal{\beta })\right] =K(\mathcal{\beta })\in
\left( 0,\infty \right) .
\end{equation*}%
In particular,%
\begin{eqnarray*}
K(\mathcal{\infty }) &=&\lim_{\min (j,n-j)\rightarrow \infty }j^{3/2}\left(
n-j\right) ^{1/2}\mathbf{E}\left[ \mathcal{V}_{j,n}(\mathcal{\infty })\right]
\\
&=&\lim_{\min (i,n-i)\rightarrow \infty }\left( n-i\right) ^{3/2}i^{1/2}%
\mathbf{E}\left[ \mathcal{H}_{i,n}(0)\right] \\
&=&\lim_{\min (i,n-i)\rightarrow \infty }\left( n-i\right) ^{3/2}i^{1/2}%
\mathbf{P}\left( \mathcal{A}_{i}(n)\right) \in (0,\infty ),
\end{eqnarray*}%
where we used the substitution $i\rightarrow n-j$ and (\ref{AsymInterm}).
Recalling (\ref{YaglomForm}) and (\ref{intro_Nu}), we conclude that
\begin{eqnarray*}
\lim_{\min (i,n-i)\rightarrow \infty }\mathbf{E}\left[ e^{-\mathcal{\beta }%
Y_{i,n}}|\mathcal{A}_{i}(n)\right] &=&1-\lim_{\min (i,n-i)\rightarrow \infty
}\frac{\mathbf{E}\left[ \mathcal{H}_{i,n}(\exp \left\{ -\mathcal{\beta }%
a_{i,n}\right\} )\right] }{\mathbf{P}\left( \mathcal{A}_{i}(n)\right) } \\
&=&1-\lim_{\min (i,n-i)\rightarrow \infty }\frac{i^{1/2}\left( n-i\right)
^{3/2}\mathbf{E}\left[ \mathcal{H}_{i,n}(\exp \left\{ -\mathcal{\beta }%
a_{i,n}\right\} )\right] }{i^{1/2}\left( n-i\right) ^{3/2}\mathbf{P}\left(
\mathcal{A}_{i}(n)\right) } \\
&=&1-\frac{K(\mathcal{\beta })}{K(\infty )}:=\Lambda \left( \mathcal{\beta }%
\right) =\mathbf{E}\left[ e^{-\beta \hat{Y}}\right] .
\end{eqnarray*}

To complete the proof of point 3) of Theorem \ref{T_conditional} it remains
to show that
\begin{equation*}
\lim_{\beta \uparrow \infty }\Lambda \left( \mathcal{\beta }\right) =\mathbf{%
P}\left( \hat{Y}=0\right) =0\text{ and }\lim_{\beta \downarrow 0}\Lambda
\left( \mathcal{\beta }\right) =\mathbf{P}\left( \hat{Y}<\infty \right) =1.
\end{equation*}%
To check the validity of these statements we first observe that, for all $%
\beta \in (0,\infty )$%
\begin{equation*}
K(\beta )=K(\infty )-K(\infty )\Lambda \left( \mathcal{\beta }\right)
=\sum_{k=0}^{\infty }C_{k}(\mathcal{\beta })+\sum_{k=1}^{\infty }J_{-k}(%
\mathcal{\beta })+\sum_{k=0}^{\infty }J_{+k}(\mathcal{\beta })\leq K(\infty
).
\end{equation*}%
Combining this estimate with (\ref{C_k_vanish}), (\ref{J_minus_vanish}) and (%
\ref{Jplus_vanish}), and applying the monotone convergence theorem we see that
\begin{equation*}
\lim_{\beta \downarrow 0}\left( K(\infty )-K(\infty )\Lambda \left( \mathcal{%
\beta }\right) \right) =K(\infty )-K(\infty )\lim_{\beta \downarrow
0}\Lambda \left( \mathcal{\beta }\right) =0.
\end{equation*}%
Thus, $\lim_{\beta \downarrow 0}\Lambda \left( \mathcal{\beta }\right) =1$
and, therefore, $\hat{Y}$ is a proper random variable.

Moreover, in view of (\ref{C_k_infinite}), (\ref{J_minus_infinity}) and (\ref%
{J_plus_infinity})
\begin{eqnarray}
&&\lim_{\beta \uparrow \infty }\left( K(\infty )-K(\infty )\Lambda \left(
\mathcal{\beta }\right) \right) =K(\infty )(1-\lim_{\beta \uparrow \infty
}\Lambda \left( \mathcal{\beta }\right) )  \notag \\
&&\qquad=\sum_{k=0}^{\infty }C_{k}(\mathcal{\infty })+\sum_{k=1}^{\infty
}J_{-k}(\mathcal{\infty })+\sum_{k=0}^{\infty }J_{+k}(\mathcal{\infty }%
)=K(\infty )=K.  \label{ExplicK}
\end{eqnarray}

These relations justify the equalities
\begin{equation*}
\mathbf{P}\left( \hat{Y}=0\right) =\lim_{\beta \uparrow \infty }\Lambda
\left( \mathcal{\beta }\right) =0.
\end{equation*}%
Thus, $\hat{Y}$ is positive with probability 1.

This completes the proof of point 3) of Theorem \ref{T_conditional}.

\section{The case of fixed $i$}

\label{sec_fixed_i}

The proof of point 2) of Theorem \ref{T_conditional} is similar to the proof
of point 3) of the theorem and is shorter.

Set, for $x\geq 1$%
\begin{equation*}
\mathcal{M}_{j}(\mathcal{\beta };x):=\frac{\bar{a}_{j}}{(1-\exp \left\{ -%
\mathcal{\beta }/\bar{a}_{j}\right\} )^{-1}+\bar{B}_{1,j+1}}\frac{\bar{b}%
_{j+1}}{\bar{b}_{j}}\frac{1}{\bar{b}_{j}+\bar{a}_{j}x}.
\end{equation*}%
We write as earlier $j$ instead of $n-i$ and taking the expectation with
respect to
\begin{equation*}
\left( \bar{X}_{j+1},\bar{X}_{j+2},...,\bar{X}_{n}\right) \overset{d}{=}%
\left( \bar{X}_{1},\bar{X}_{2},...,\bar{X}_{i}\right)
\end{equation*}%
and recalling (\ref{DefNu}) obtain that
\begin{equation*}
\mathbf{E}\left[ \mathcal{H}_{i,n}(\exp \left\{ -\mathcal{\beta }%
a_{i,n}\right\} )\right] =\mathbf{E}\left[ \mathcal{V}_{j,n}(\mathcal{\beta }%
)\right] =\mathbf{E}\left[ \Upsilon _{j}(\mathcal{\beta };\bar{b}_{i})\right]
,
\end{equation*}%
where%
\begin{equation*}
\Upsilon _{j}(\mathcal{\beta };x):=\mathbf{E}\left[ \mathcal{M}_{j}(\mathcal{%
\beta };x)\right] .
\end{equation*}%
Note that in view of (\ref{AboveRW}) for all $\beta \geq 0$ and $x\geq 1$%
\begin{equation}
\mathbf{E}\left[ \mathcal{M}_{j}(\mathcal{\beta };x)\right] \leq \mathbf{E}%
\left[ \frac{\bar{a}_{j}}{\bar{b}_{j}}\frac{1}{\bar{b}_{j}+\bar{a}_{j}}%
\right] \leq \mathbf{E}\left[ e^{-\bar{S}_{j}}e^{\bar{S}_{\tau (j-1)}}e^{%
\bar{S}_{\tau (j)}}\right] \leq \frac{C}{j^{3/2}}.  \label{AboveRW_j}
\end{equation}%
Recall the definition of $\mathcal{K}_1$ and $\mathcal{K}_2$ in (\ref%
{defK1K2}). Now we fix some positive integer $N$ and use the decomposition
\begin{eqnarray*}
\Upsilon _{j}(\mathcal{\beta };x) &=&\mathbf{E}\left[ \mathcal{M}_{j}(%
\mathcal{\beta };x);\bar{\tau}(j)\in \mathcal{K}_{1}\right] +\mathbf{E}\left[
\mathcal{M}_{j}(\mathcal{\beta };x);\bar{\tau}(j)<N\right]  \notag \\
&+&\mathbf{E}\left[ \mathcal{M}_{j}(\mathcal{\beta };x);\bar{\tau}(j)\in
(j-N,j]\right] .  \label{Decomposition2}
\end{eqnarray*}%

First we observe that in view of (\ref{AsymConditional}) for any $%
\varepsilon >0$ one can find $N_{0}=N_{0}(\varepsilon )$ such that
\begin{eqnarray}
&&\mathbf{E}\left[ \mathcal{M}_{j}(\mathcal{\beta };x);\bar{\tau}(j)\in
\mathcal{K}_{1}\right] \leq \mathbf{E}\left[e^{-\bar{S}_{j}};\bar{\tau}%
(j)\in \mathcal{K}_{1}\right]  \notag \\
&&\qquad\qquad\quad=\sum_{r=N}^{j-N}\mathbf{E}\left[ e^{\bar{S}_{r}-\bar{S}_{j}}e^{-\bar{S}%
_{r}};\bar{\tau}(j)=r\right]  \notag \\
&&\qquad\qquad\quad=\sum_{r=N}^{j-N}\mathbf{E}\left[ e^{-\bar{S}_{r}};\bar{\tau}(r)=r\right]
\mathbf{E}\left[ e^{-\bar{S}_{j-r}};L_{j-r}\geq 0\right]  \notag \\
&&\qquad\qquad\quad\leq\frac{C}{j^{3/2}}\sum_{r=N}^{\infty }\frac{1}{r^{3/2}}\leq \frac{%
\varepsilon }{j^{3/2}}  \label{FirstMnegl}
\end{eqnarray}%
for all sufficiently large $j$ and $N\geq N_{0}$.

Now we evaluate $\mathbf{E}\left[ \mathcal{M}_{j}(\mathcal{\beta };x);\bar{%
\tau}(j)<N\right] $. By conditioning on the trajectory of the random walk $%
\mathbf{\bar{S}}$ until time $k$ we obtain
\begin{equation*}
\mathbf{E}\left[ \mathcal{M}_{j}(\mathcal{\beta };x);\bar{\tau}(j)=k\right] =%
\mathbf{E}\left[ e^{-\bar{S}_{k}}\Psi _{j-k}\left( \mathcal{\beta };x,e^{-%
\bar{S}_{k}},\bar{b}_{k}\right) ;\bar{\tau}(k)=k\right]
\end{equation*}%
where%
\begin{eqnarray*}
&&\Psi _{j}\left( \mathcal{\beta };x,u,q\right) \\
&:&=\mathbf{E}\left[ e^{-\bar{S}_{j}}\frac{1}{(1-\exp \left\{ -\mathcal{%
\beta }/u\bar{a}_{j}\right\} )^{-1}+q-1+u\bar{b}_{j+1}}\frac{q+u\bar{b}_{j+1}%
}{q+u\bar{b}_{j}}\frac{1}{q+u(\bar{b}_{j}+\bar{a}_{j}x)};\bar{L}_{j}\geq 0%
\right] .
\end{eqnarray*}%
Recall that $t=\left[ j/2\right] $, as well as the definition of $\tilde{H}$
in (\ref{NewT}). A direct computation gives that
\begin{equation*}
\Psi _{j}\left( \mathcal{\beta };x,u,q\right) =\mathbf{E}\left[ e^{-\bar{S}%
_{j}}\varphi _{\mathcal{\beta },x,u,q}(G_{t},\tilde{%
H}_{t+1,j-1},S_{j});\bar{L}_{j}\geq 0\right] ,
\end{equation*}%
where, for $\beta \in (0,\infty )$
\begin{eqnarray*}
\varphi _{\mathcal{\beta },x,u,q}(g,h,s) &:&=\frac{1}{(1-\exp \left\{ -%
\mathcal{\beta }/ue^{-s}\right\} )^{-1}+q-1+u\left( g+e^{-s}(h+1)\right) } \\
&&\times \frac{q+u\left( g+e^{-s}(h+1)\right) }{q+u\left( g+e^{-s}h\right) }%
\frac{1}{q+u\left( g+e^{-s}(h+x)\right) },
\end{eqnarray*}%
and, for $\beta =\infty $%
\begin{equation}
\varphi _{\mathcal{\infty },x,u,q}(g,h,s):=\frac{1}{q+u\left(
g+e^{-s}h\right) }\frac{1}{q+u\left( g+e^{-s}(h+x)\right) }\leq 1
\label{BoundPhi}
\end{equation}%
if $q\geq 1$.

Repeating now the arguments similar to those used to prove (\ref{Psi_j,n})
and applying Theorem 2.7 from \cite{ABKV} we obtain
\begin{equation*}
\lim_{j\rightarrow \infty }\frac{\Psi _{j}\left( \mathcal{\beta }%
;x,u,q\right) }{\mathbf{E}[e^{-\bar{S}_{j}};\bar{L}_{j}\geq 0]\ }=\Psi
_{\infty }\left( \mathcal{\beta };x,u,q\right) ,  \label{Psi_j_infty}
\end{equation*}%
where for $\beta \in (0,\infty ]$%
\begin{equation*}
\Psi _{\infty }\left( \mathcal{\beta };x,u,q\right) :=\int^{(3)}\varphi
_{\beta ,x,u,q}(g,h,-s)\mathbf{\bar{P}}^{+}\left( G_{\infty }\in dg\right)
\mathbf{\bar{P}}_{s}^{-}\left( H_{1,\infty }\in dh\right) \bar{\nu}_{1}(ds).
\end{equation*}%
Note that in view of (\ref{BoundPhi})%
\begin{equation*}
\Psi _{\infty }\left( \mathcal{\beta };x,u,q\right) \leq \int^{(3)}\mathbf{%
\bar{P}}^{+}\left( G_{\infty }\in dg\right) \mathbf{\bar{P}}_{s}^{-}\left(
H_{1,\infty }\in dh\right) \bar{\nu}_{1}(ds)\leq \int \bar{\nu}_{1}(ds)=1.
\end{equation*}

Since \
\begin{equation*}
\lim_{\beta \downarrow 0}\varphi _{\beta ,x,u,q}(g,h,-s)=0
\end{equation*}%
and $\varphi _{\beta ,x,u,q}$ is monotone decreasing as $\beta \downarrow 0$%
, we conclude by the dominated convergence theorem that%
\begin{equation*}
\lim_{\beta \downarrow 0}\Psi _{\infty }\left( \mathcal{\beta };x,u,q\right)
=0.  \label{Psi_zero2}
\end{equation*}

Furthermore, setting
\begin{equation*}
C_{k}^{\prime }(\mathcal{\beta },x):=\mathbf{E}\left[ e^{-\bar{S}_{k}}\Psi
_{\infty }\left( \mathcal{\beta };x,e^{-\bar{S}_{k}},\bar{B}_{1,k}\right) ;%
\bar{\tau}(k)=k\right] ,\quad \beta \in (0,\infty ],  \label{DefC_kx}
\end{equation*}%
we have, again by monotonicity of $\Psi _{\infty }\left( \mathcal{\beta }%
;\cdot ,\cdot ,\cdot \right) $ in $\beta $ and the dominated convergence
theorem that
\begin{equation}
\lim_{\beta \downarrow 0}C_{k}^{\prime }(\mathcal{\beta },x)=0
\label{C_kx_Vanish}
\end{equation}%
and
\begin{eqnarray*}
\lim_{\beta \uparrow \infty }C_{k}^{\prime }(\mathcal{\beta },x) &=&\mathbf{E%
}\left[ e^{-\bar{S}_{k}}\Psi _{\infty }\left( \mathcal{\infty };x,e^{-\bar{S}%
_{k}},\bar{B}_{1,k}\right) ;\bar{\tau}(k)=k\right]  \notag \\
&=&C_{k}^{\prime }(\infty ,x)<\infty .  \label{C_kx_infinity}
\end{eqnarray*}%
Applying the dominated convergence theorem once more we conclude by (\ref%
{AsymConditional}) that for any fixed $k$, as $j\rightarrow \infty $%
\begin{eqnarray}
&&\mathbf{E}\left[ e^{-\bar{S}_{k}}\Psi _{j-k}\left( \mathcal{\beta };x,e^{-%
\bar{S}_{k}},\bar{B}_{1,k}\right) ;\bar{\tau}(k)=k\right]  \notag \\
&&\qquad \qquad \sim \mathbf{E}\left[ e^{-\bar{S}_{k}}\Psi _{\infty }\left(
\mathcal{\beta };x,e^{-\bar{S}_{k}},\bar{B}_{1,k}\right) ;\bar{\tau}(k)=k%
\right] \mathbf{E}[e^{-\bar{S}_{j}};\bar{L}_{j}\geq 0]  \notag \\
&&\qquad \qquad \qquad \qquad \sim \frac{K_{4}C_{k}^{\prime }(\mathcal{\beta
},x)}{j^{3/2}}.  \label{TermLeft2}
\end{eqnarray}%
Summing (\ref{TermLeft2}) over the $k$'s in $[0,N-1]$, gives, as $%
j\rightarrow \infty $
\begin{equation}
j^{3/2}\mathbf{E}\left[ \mathcal{M}_{j}(\mathcal{\beta };x);\bar{\tau}(j)<N%
\right] \sim K_{4}C^{\prime }(\mathcal{\beta },x,N)  \label{LeftPart4}
\end{equation}%
where%
\begin{equation*}
C^{\prime }(\mathcal{\beta },x,N):=\sum_{k=0}^{N-1}\mathbf{E}\left[ e^{-\bar{S}%
_{k}}\Psi _{\infty }\left( \mathcal{\beta };x,e^{-\bar{S}_{k}},\bar{B}%
_{1,k}\right) ;\bar{\tau}(k)=k\right] .  \label{LeftPart5}
\end{equation*}%
We now evaluate $\mathbf{E}\left[ \mathcal{M}_{j}(\mathcal{\beta };x);\bar{%
\tau}(n)\in (j-N,j]\right] $ and write
\begin{eqnarray*}
&&\mathbf{E}\left[ \mathcal{M}_{j}(\mathcal{\beta };x);\bar{\tau}(n)=j-k%
\right] \\
&=&\mathbf{E}\left[ e^{\bar{S}_{j-k}}\varphi _{\beta ,x}^{-}(G_{t},\tilde{H}_{t,j-k-1,1},\bar{S}_{j-k},\bar{S}_{j-k}-\bar{S}%
_{j},\hat{T}_{j-k,j},\hat{T}_{j-k,j-1});\bar{\tau}(n)=j-k\right] ,
\end{eqnarray*}

where, for $\beta \in (0,\infty )$%
\begin{eqnarray*}
\varphi _{\beta ,x}^{-}(g,h,s_{1},s_{2},t_{1},t_{2}) &:=&\frac{e^{s_{2}}}{%
e^{s_{1}}(1-\exp \left\{ -\mathcal{\beta }e^{-s_{2}}e^{s_{1}}\right\}
)^{-1}+e^{s_{1}}(g-1)+h+t_{1}} \\
&&\times \frac{e^{s_{1}}g+h+t_{1}}{e^{s_{1}}g+h+t_{2}}\frac{1}{%
e^{s_{1}}g+h+t_{2}+x}
\end{eqnarray*}%
and, by monotonicity
\begin{eqnarray}
\varphi _{\beta ,x}^{-}(g,h,s_{1},s_{2},t_{1},t_{2}) &\leq &\varphi _{\infty
,x}^{-}(g,h,s_{1},s_{2},t_{1},t_{2})  \notag \\
&:&=\frac{e^{s_{2}}}{e^{s_{1}}g+h+t_{2}}\frac{1}{e^{s_{1}}g+h+t_{2}+x}\leq
1  \label{Phi_minus_bound_x}
\end{eqnarray}%
in the domain%
\begin{equation*}
\left\{ g\geq 1,h\geq 0,s_{1}\leq 0,s_{2}\leq 0,t_{1}\geq 0,t_{2}\geq
0\right\} \cap \left\{ e^{s_{2}}\leq
e^{s_{1}}g+h+t_{2},e^{s_{1}}g+h+t_{2}+x\geq 1\right\} .
\end{equation*}%
Since
\begin{eqnarray*}
G_{t} &\rightarrow &G_{\infty }\quad \mathbf{\bar{P}}%
_{x}^{+}-a.s.\,\forall x\geq 0,\,(H_{1,j},S_{k})\rightarrow (H_{1,\infty
},S_{k})\quad \mathbf{\bar{P}}^{-}-a.s.,\, \\
(T_{k-1},T_{k}) &\rightarrow &(T_{k-1},T_{k}),\quad \mathbf{\bar{P}}^{+}-a.s.
\end{eqnarray*}%
as $j\rightarrow \infty $ we may apply Lemma \ref{L_pr6New} in
\cite{ABKV} and conclude that, for each $k$, as $j\rightarrow \infty $
\begin{eqnarray}
\mathbf{E}\left[ \mathcal{M}_{j}(\mathcal{\beta },x);\bar{\tau}(j)=j-k\right]
&\sim &J_{-k}(\mathcal{\beta },x)\mathbf{E}[e^{\bar{S}_{j-k}};\bar{\tau}%
(j-k)=j-k]  \notag \\
&\sim &\frac{K_{3}J_{-k}(\mathcal{\beta },x)}{j^{3/2}},  \label{Term2a_x}
\end{eqnarray}%
where
\begin{eqnarray*}
J_{-k}(\mathcal{\beta },x) &=&\mathbf{P}\left( \bar{L}_{k}\geq 0\right)
\int^{(6)}\varphi _{\beta ,x}^{-}(g,h,-s_{1},-s_{2},t_{1},t_{2})  \notag \\
&&\times \mathbf{\bar{P}}_{s_{1}}^{+}\left( G_{\infty }\in dg\right) \mathbf{%
\bar{P}}^{+}\left( (T_{k-1},T_{k})\in d(t_{1},t_{2})\right)  \notag \\
&&\times \mathbf{\bar{P}}^{-}\left( (H_{1,\infty },S_{k})\in
d(h,s_{2})\right) \bar{\mu}_{1}(ds_{1}).  \label{Jminus_x}
\end{eqnarray*}

Monotonicity of $\varphi _{\beta
,x}^{-}(g,h,-s_{1},-s_{2},t_{1},t_{2}) $ in $\beta $ and the dominated
convergence theorem show that, for each fixed $k$%
\begin{equation}
\lim_{\beta \downarrow 0}J_{-k}(\mathcal{\beta },x)=0
\label{J_minus_vanish_x}
\end{equation}%
and%
\begin{equation*}
\lim_{\beta \uparrow \infty }J_{-k}(\mathcal{\beta },x)=J_{-k}(\mathcal{%
\infty },x)<\infty .  \label{J_minus_infinity_x}
\end{equation*}%
Summing (\ref{Term2a_x}) over the $k$'s in $[0,N-1]$, gives, as $%
j\rightarrow \infty $
\begin{equation}
j^{3/2}\mathbf{E}\left[ \mathcal{M}_{j}(\mathcal{\beta },x);\bar{\tau}(j)\in
(j-N,j]\right] \sim K_{3} J_{-}(\mathcal{\beta },x,N)= K_{3}
\sum_{k=0}^{N-1}J_{-k}(\mathcal{\beta },x).  \label{LeftTerm2_x}
\end{equation}%
Combining (\ref{LeftPart4}) and (\ref{LeftTerm2_x}) with (\ref{AboveRW_j})
and (\ref{FirstMnegl}) shows that
\begin{equation*}
\lim_{j\rightarrow \infty }j^{3/2}\Upsilon _{j}(\mathcal{\beta }%
;x)=\lim_{j\rightarrow \infty }j^{3/2}\mathbf{E}\left[ \mathcal{M}_{j}(%
\mathcal{\beta },x)\right] =K(\mathcal{\beta },x)
\end{equation*}%
where, for all $\beta \in (0,\infty ]$%
\begin{equation*}
K(\mathcal{\beta },x):=K_{4}\sum_{k=0}^{\infty }C_{k}^{\prime }(\mathcal{%
\beta },x)+K_{3}\sum_{k=0}^{\infty }J_{-k}(\mathcal{\beta },x)\leq K(%
\mathcal{\infty },x)\in \left( 0,\infty \right) .
\end{equation*}%
Hence, using the dominated convergence theorem we conclude that
\begin{eqnarray*}
\lim_{j\rightarrow \infty }j^{3/2}\mathbf{E}\left[ \mathcal{H}_{i,n}(\exp
\left\{ -\mathcal{\beta }a_{i,n}\right\} )\right] &=&\lim_{j\rightarrow
\infty }j^{3/2}\mathbf{E}\left[ \Upsilon _{j}(\mathcal{\beta };\bar{b}_{i})%
\right] \\
&=&\mathbf{E}\left[ \lim_{j\rightarrow \infty }j^{3/2}\Upsilon _{j}(\mathcal{%
\beta };\bar{b}_{i})\right] =\mathbf{E}\left[ K(\mathcal{\beta };\bar{b}_{i})%
\right] .
\end{eqnarray*}

Observe that in view of (\ref{Phi_minus_bound_x}) the
arguments above are valid for the case $\beta =\infty,$ i.e. for the case
when the expression $1-\exp \left\{ -\beta u\bar{a}_{j}\right\} $ is
everywhere replaced by 1. This corresponds to studying the asymptotic
behavior of \ $\mathbf{P}\left( \mathcal{A}_{i}(n)\right) $ as $n\rightarrow
\infty $. Therefore, $w_{i}=\mathbf{E}\left[ K(\mathcal{\infty };\bar{b}_{i})%
\right] $ in (\ref{AsumBegin}). Thus, for each fixed $i$ and $j=n-i$%
\begin{eqnarray*}
\lim_{n\rightarrow \infty }\mathbf{E}\left[ e^{-\mathcal{\beta }Y_{i,n}}|%
\mathcal{A}_{i}(n)\right] &=&1-\lim_{n\rightarrow \infty }\frac{\mathbf{E}%
\left[ \mathcal{H}_{i,n}(\exp \left\{ -\mathcal{\beta }a_{i,n}\right\} )%
\right] }{\mathbf{P}\left( \mathcal{A}_{i}(n)\right) }  \notag \\
&=&1-\frac{\lim_{n\rightarrow \infty }(n-i)^{3/2}\mathbf{E}\left[ \Upsilon
_{j}(\mathcal{\beta };\bar{b}_{i})\right] }{\lim_{n\rightarrow \infty
}\left( n-i\right) ^{3/2}\mathbf{P}\left( \mathcal{A}_{i}(n)\right) }  \notag
\\
&=&1-\frac{\mathbf{E}\left[ K(\mathcal{\beta };\bar{b}_{i})\right] }{\mathbf{%
E}\left[ K(\mathcal{\infty };\bar{b}_{i})\right] }:=\Lambda _{i}\left(
\mathcal{\beta }\right) =\mathbf{E}\left[ e^{-\beta \hat{Y}_{i,\infty }}%
\right] .  \label{Zero_j}
\end{eqnarray*}%
Again by monotonicity of $K(\mathcal{\beta },x)$ and (\ref{C_kx_Vanish}) and
(\ref{J_minus_vanish_x}) it follows that%
\begin{equation*}
\lim_{\beta \downarrow 0}K(\mathcal{\beta },x)=0
\end{equation*}%
and%
\begin{equation*}
\lim_{\beta \downarrow 0}\mathbf{E}\left[ K(\mathcal{\beta };\bar{b}_{i})%
\right] =\mathbf{E}\left[ \lim_{\beta \downarrow 0}K(\mathcal{\beta };\bar{b}%
_{i})\right] =0.
\end{equation*}%
Thus,%
\begin{equation*}
\mathbf{P}\left( \hat{Y}_{i,\infty }<\infty \right) =1.
\end{equation*}%
On the other hand, the inequality $K(\mathcal{\beta };x)\leq K(\mathcal{%
\infty };x)\leq K(\mathcal{\infty };1)<\infty $ and the dominated
convergence theorem give%
\begin{equation*}
\lim_{\beta \uparrow \infty }\mathbf{E}\left[ K(\mathcal{\beta };\bar{b}_{i})%
\right] =\mathbf{E}\left[ \lim_{\beta \uparrow \infty }K(\mathcal{\beta };%
\bar{b}_{i})\right] =\mathbf{E}\left[ K(\mathcal{\infty };\bar{b}_{i})\right]
\end{equation*}%
implying
\begin{equation*}
\lim_{\beta \uparrow \infty }\Lambda \left( \beta \right) =\mathbf{P}\left(
\hat{Y}_{i,\infty }=0\right) =0.
\end{equation*}%
This completes the proof of point 2) of Theorem \ref{T_conditional}.\\

The authors thank the reviewer, whose comments allow to eliminate a number of inaccuracies contained in the original version of the article.

\end{document}